\documentclass[leqno,12pt]{article}%
\usepackage{amsfonts}
\usepackage{amsmath}
\usepackage{amssymb}
\usepackage{graphicx}%
\newtheorem{theorem}{Theorem}
\newtheorem{acknowledgement}[theorem]{Acknowledgement}

\newtheorem{definition}[theorem]{Definition}

\newtheorem{notation}[theorem]{Notation}

\newtheorem{proposition}[theorem]{Proposition}
\newtheorem{remark}[theorem]{Remark}

\newenvironment{proof}[1][Proof]{\noindent\textbf{#1.} }{\ \rule{0.5em}{0.5em}}
\begin{document}

\title{{\huge The Fitzpatrick function - a bridge between convex analysis and
multivalued stochastic differential equations}\thanks{{\small The work for
this paper was supported by founds from the Grant CNCSIS nr. 1373/2007 and
Grant CNCSIS nr. 1156/2005.}}}
\author{{\large Aurel R\u{a}\c{s}canu}$^{~1}\quad\quad${\large Eduard Rotenstein}%
$^{~2}$}
\date{}
\maketitle

\begin{abstract}
Using the Fitzpatrick function, we characterize the solutions for different
classes of deterministic and stochastic differential equations driven by
maximal monotone operators (or in particular subdifferential operators) as the
minimum point of a suitably chosen convex lower semicontinuous function. Such
technique provides a new approach for the existence of the solutions for the
considered equations.\bigskip

\end{abstract}

\medskip

\footnotetext[1]{{\small Department of Mathematics, "Al. I. Cuza" University,
Bd. Carol no.9-11 \& "Octav Mayer" Mathematics Institute of the Romanian
Academy, Bd. Carol I, no.8, Romania, \newline e-mail: aurel.rascanu@uaic.ro}}

\footnotetext[2]{{\small Department of Mathematics, "Al. I. Cuza" University,
Bd. Carol no.9-11, Ia\c{s}i, Rom\^{a}nia, e-mail: eduard.rotenstein@uaic.ro}}


\section{Preliminaries. Notations}

The Fitzpatrik function proved to be a very useful tool of the convex analysis
in the study of maximal monotone operators. In our paper this function is used
for deterministic and stochastic differential equations driven by multivalued
maximal monotone operators. We will show how we can reduce the existence
problem for stochastic differential equations of the following types:\bigskip

\begin{itemize}
\item \textit{forward case}%
\begin{equation}
\left\{
\begin{array}
[c]{l}%
dX_{t}+A\left(  X_{t}\right)  \left(  dt\right)  \ni F\left(  t,X_{t}\right)
dt+G\left(  t,X_{t}\right)  dW_{t}~,\smallskip\\
X_{0}=\xi,\ t\in\left[  0,T\right]  \quad\quad\text{and}%
\end{array}
\right.  \label{fc1}%
\end{equation}

\item \textit{backward case}%
\begin{equation}
\left\{
\begin{array}
[c]{l}%
-dY_{t}+A\left(  Y_{t}\right)  dt\ni H\left(  t,Y_{t},Z_{t}\right)
dt-Z_{t}dW_{t}~,\smallskip\\
Y_{T}=\xi,\ t\in\left[  0,T\right]
\end{array}
\right.  \label{bc1}%
\end{equation}

\end{itemize}

\noindent to a minimizing problem for convex lower semicontinuous functions.

Usually, existence results are obtained via a penalized problem with Yosida's
approximation operator $A_{\varepsilon}:=[I-\left(  I+\varepsilon A\right)
^{-1}]/\varepsilon$.

For the forward equation (\ref{fc1}), by studying first a generalized Skorohod
problem%
\[
\left\{
\begin{array}
[c]{l}%
dx\left(  t\right)  +A\left(  x\left(  t\right)  \right)  \left(  dt\right)
\ni f\left(  t\right)  dt+dm\left(  t\right)  ,\smallskip\\
x\left(  0\right)  =x_{0},\ t\in\left[  0,T\right]  .
\end{array}
\right.
\]
the existence of the solution is obtained (see Bensoussan \& R\u{a}\c{s}canu
\cite{BR}, R\u{a}\c{s}canu \cite{R}, or Asiminoaei \& R\u{a}\c{s}canu
\cite{AR}) in the general case of a maximal monotone operator.

For backward stochastic differential equations the existence problem (see
Pardoux \& Rascanu \cite{PR}) is solved only in the case of $A=\partial
\varphi$ (the subdifferential of a lower semicontinuous convex function) and
it is an open problem in the general case. That is \textit{the reason and the
main motivation} to find an approach via convex analysis.

In 1988, in the paper \cite{F88}, Fitzpatrick proved that any maximal monotone
operator can be represented by a convex function; he explicitly defined the
minimal convex representation. The connection between maximal monotone
operators and convex functions was also approached 13 years later by
Martinez-Legaz \& Thera in \cite{MT}, Burachik \& Svaiter in \cite{BS} and
Burachik \& Fitzpatrick in \cite{BF}. Since these last three papers,
Fitzpatrick's results have been the subject of intense research (J.P.
Revalski, M. Thera, R.S. Burachik, B.F. Svaiter, J.-P. Penot, S. Simons, C.
Z\u{a}linescu, J.-E. Martinez-Legaz etc.). Their results stay in the domain of
nonlinear operators: properties, characterizations, new classes of monotone operators.

Using the idea of Fitzpatrick function we can reduce the existence problems
for stochastic equations of the form (\ref{fc1}) or (\ref{bc1}) to a
minimizing problem of a convex lower semicontinuous function. Inspired by the
studies of Gy\"{o}ngy \& Mart\'{\i}nez \cite{GM}, we present a new approach
for solving the existence problem for stochastic differential equations with
maximal monotone operator. In this paper we will identify the solutions of
different types of forward and backward multivalued stochastic differential
equations with the minimum points of a suitably chosen convex lower
semicontinuous functionals.

The paper is organized as follows. In the first section we present some basic
properties of the Fitzpatrick's function and we will introduce the stochastic
framework that will be used. The next section contains a Fitzpatrick function
approach for the study of a generalized Skorohod problem as well of forward
and backward stochastic differential equations, while Section 3 is dedicated
to the case of forward and backward stochastic variational inequalities.

\subsection{On Fitzpatrick's function}

Let $\left(  \mathbb{X},\left\Vert .\right\Vert \right)  $ be a real Banach
space and $\left(  \mathbb{X}^{\ast},\left\Vert .\right\Vert _{\ast}\right)  $
be its dual. For $x^{\ast}\in\mathbb{X}^{\ast}$ and $x\in\mathbb{X}$ we denote
$x^{\ast}\left(  x\right)  $ (the value of $x^{\ast}$ in $x$) by $\left\langle
x,x^{\ast}\right\rangle $ or $\left\langle x^{\ast},x\right\rangle .$

If $A:\mathbb{X}\rightrightarrows\mathbb{X}^{\ast}$ is a point-to-set operator
(from $\mathbb{X}$ to the family of subsets of $\mathbb{X}^{\ast}$), then
$Dom\left(  A\right)  :=\left\{  x\in\mathbb{X}:A\left(  x\right)
\neq\emptyset\right\}  $ and $R\left(  A\right)  =\left\{  x^{\ast}:\;\exists
x\in Dom\left(  A\right)  \text{ s.t. }x^{\ast}\in A\left(  x\right)
\right\}  .$ We shall always assume that the operator $A$ is proper, i.e.
$Dom\left(  A\right)  \neq\emptyset.$ Usually the operator $A$ is identified
with its graph $gr\left(  A\right)  =\{\left(  x,x^{\ast}\right)
\in\mathbb{X\times X}^{\ast}:x\in Dom\left(  A\right)  ,\;x^{\ast}\in A\left(
x\right)  \}.$

The operator $A:\mathbb{X}\rightrightarrows\mathbb{X}^{\ast}$ is a monotone
operator ($A\subset\mathbb{X\times X}^{\ast}$ is a monotone set) if%
\[
\left\langle x-y,x^{\ast}-y^{\ast}\right\rangle \geq0,\ \forall\left(
x,x^{\ast}\right)  ,\left(  y,y^{\ast}\right)  \in A.
\]
A monotone operator (set) is maximal monotone if it is not properly contained
in any other monotone operator (set). Clearly if $A$ is maximal monotone and
$\left(  y,y^{\ast}\right)  \in\mathbb{X}\times\mathbb{X}^{\ast}$ then%
\[
\inf\limits_{\left(  u,u^{\ast}\right)  \in A}\left\langle y-u,y^{\ast
}-u^{\ast}\right\rangle \geq0\quad\Longleftrightarrow\quad\left(  y,y^{\ast
}\right)  \in A.
\]

Given a function $\psi:\mathbb{X\rightarrow]-}\infty,+\infty]$ we denote
$Dom\left(  \psi\right)  :=\left\{  x\in\mathbb{X}:\psi\left(  x\right)
<\infty\right\}  .$ We say that $\psi$ is proper if $Dom\left(  \psi\right)
\neq\emptyset.$ The subdifferential $\partial\psi:\mathbb{X}\rightrightarrows
\mathbb{X}^{\ast}$ is defined by%
\[
\left(  x,x^{\ast}\right)  \in\partial\psi\quad\text{if}\quad\left\langle
y-x,x^{\ast}\right\rangle +\psi\left(  x\right)  \leq\psi\left(  y\right)
,\ \forall y\in\mathbb{X}.
\]
It is well known that: if $\psi$ is a proper convex l.s.c. function, then
$\partial\psi:\mathbb{X}\rightrightarrows\mathbb{X}^{\ast}$ is a maximal
monotone operator.

Let $\psi:\mathbb{X\rightarrow]-}\infty,+\infty]$ be a proper function. The
conjugate of $\psi$ is the function\newline$\psi^{\ast}:\mathbb{X}^{\ast
}\mathbb{\rightarrow]-}\infty,+\infty],$
\[
\psi^{\ast}\left(  x^{\ast}\right)  :=\sup\left\{  \left\langle u,x^{\ast
}\right\rangle -\psi\left(  u\right)  :u\in\mathbb{X}\right\}  .
\]

Remark that, if $h:\mathbb{X\times X}^{\ast}\mathbb{\rightarrow]-}%
\infty,+\infty]$, then $h^{\ast}:\mathbb{X}^{\ast}\mathbb{\times X}^{\ast\ast
}\mathbb{\rightarrow]-}\infty,+\infty]$ and, for any\newline$\left(  x^{\ast
},x\right)  \in\mathbb{X}^{\ast}\mathbb{\times X}$, $h^{\ast}\left(  x^{\ast
},x\right)  $ is well defined by identifying $\mathbb{X}$ with its image under
canonical injection of $\mathbb{X}$ into $\mathbb{X}^{\ast\ast}$, that is,
every $x\in\mathbb{X}$ can be seen as a function $x:\mathbb{X}^{\ast
}\rightarrow\mathbb{R}$ defined by $x(x^{\ast})=x^{\ast}(x)=\left\langle
x,x^{\ast}\right\rangle $. For a complete study on maximal monotone operators,
one can consult Barbu \cite{B75} or Br\'{e}zis \cite{B}.

\begin{definition}
Given a monotone operator $A:\mathbb{X}\rightrightarrows\mathbb{X}^{\ast}$,
the associated Fitzpatrick function is defined as $\mathcal{H}=\mathcal{H}%
_{A}:\mathbb{X}\times\mathbb{X}^{\ast}\rightarrow\mathbb{]-}\infty,+\infty]$,%
\begin{equation}%
\begin{array}
[c]{ccc}%
\mathcal{H}\left(  x,x^{\ast}\right)  & := & \left\langle x,x^{\ast
}\right\rangle -\inf\left\{  \left\langle x-u,x^{\ast}-u^{\ast}\right\rangle
:\left(  u,u^{\ast}\right)  \in A\right\}  \smallskip\\
& = & \sup\left\{  \left\langle u,x^{\ast}\right\rangle +\left\langle
x,u^{\ast}\right\rangle -\left\langle u,u^{\ast}\right\rangle :\left(
u,u^{\ast}\right)  \in A\right\}
\end{array}
\label{F-def}%
\end{equation}

\end{definition}

\noindent Clearly $\mathcal{H}\left(  x,x^{\ast}\right)  \leq\left\langle
x,x^{\ast}\right\rangle $,$\ $for all $(x,x^{\ast})\in A$ and, as supremum of
convex strongly (and $(w,w^{\ast})$) continuous functions, $\mathcal{H}%
=\mathcal{H}_{A}:\mathbb{X}\times\mathbb{X}^{\ast}\rightarrow\mathbb{]-}%
\infty,+\infty]$ is a proper convex strongly (and $(w,w^{\ast})$) l.s.c.
function. Usually, we shall consider on $\mathbb{X}$ the strong topology and,
on $\mathbb{X}^{\ast}$ the $w^{\ast}$-topology; in this case, $\mathcal{H}$ is
also a l.s.c. function. Whenever is necessary, we will consider the
Fitzpatrick function $\mathcal{H}$ restricted at $\mathbb{U}\times\mathbb{V}$,
with $\mathbb{U}\subset\mathbb{X}$ and $\mathbb{V\subset X}^{\ast}$.

Let $\left(  x^{\ast},x\right)  \in\partial\mathcal{H}\left(  u,u^{\ast
}\right)  $. Then, from the definition of a subdifferential operator, we have%
\[
\left\langle \left(  x^{\ast},x\right)  ,\left(  z,z^{\ast}\right)  -\left(
u,u^{\ast}\right)  \right\rangle +\mathcal{H}\left(  u,u^{\ast}\right)
\leq\mathcal{H}\left(  z,z^{\ast}\right)  ,\ \forall\left(  z,z^{\ast}\right)
\in\mathbb{X}^{\ast\ast}\times\mathbb{X}^{\ast},
\]
or, equivalently,%
\begin{equation}%
\begin{array}
[c]{l}%
\left\langle u-x,u^{\ast}-x^{\ast}\right\rangle -\inf\left\{  \left\langle
u-y,u^{\ast}-y^{\ast}\right\rangle :\left(  y,y^{\ast}\right)  \in A\right\}
\medskip\\
\quad\leq\left\langle z-x,z^{\ast}-x^{\ast}\right\rangle -\inf\left\{
\left\langle z-y,z^{\ast}-y^{\ast}\right\rangle :\left(  y,y^{\ast}\right)
\in A\right\}  ,\ \forall\left(  z,z^{\ast}\right)  \in\mathbb{X}^{\ast\ast
}\times\mathbb{X}^{\ast}.
\end{array}
\label{F3}%
\end{equation}
Since the operator $A$ is a maximal monotone one, then%
\[%
\begin{array}
[c]{l}%
\inf\left\{  \left\langle u-y,u^{\ast}-y^{\ast}\right\rangle :\left(
y,y^{\ast}\right)  \in A\right\}  \leq0\quad\text{and}\smallskip\\
\inf\left\{  \left\langle z-y,z^{\ast}-y^{\ast}\right\rangle :\left(
y,y^{\ast}\right)  \in A\right\}  =0,\ \forall\left(  z,z^{\ast}\right)  \in
A\text{;}%
\end{array}
\]
consequently, we have%
\begin{equation}
\left(  x^{\ast},x\right)  \in\partial\mathcal{H}\left(  u,u^{\ast}\right)
\ \Longrightarrow\ \left\langle u-x,u^{\ast}-x^{\ast}\right\rangle \leq
\inf\left\{  \left\langle z-x,z^{\ast}-x^{\ast}\right\rangle :\left(
z,z^{\ast}\right)  \in A\right\}  . \label{F2}%
\end{equation}
Also, by the monotonicity of $A$, from (\ref{F3}) follows%
\[
\left(  x,x^{\ast}\right)  \in A\ \Longrightarrow\ \left(  x^{\ast},x\right)
\in\partial\mathcal{H}\left(  x,x^{\ast}\right)  .
\]

\noindent Hence, if $A:\mathbb{X}\rightrightarrows\mathbb{X}^{\ast}$ is a
maximal monotone operator, then $\mathcal{H}_{A}$ characterizes $A$ as follows.

\begin{theorem}
[Fitzpatrick](see Fitzpatrick \cite{F88}, Simons \& Z\u{a}linescu \cite{SZ})
Let $A:\mathbb{X}\rightrightarrows\mathbb{X}^{\ast}$ be a maximal monotone
operator and $\mathcal{H}$ its associated Fitzpatrick function. Then, for all
$(x,x^{\ast})\in\mathbb{X}\times\mathbb{X}^{\ast}$,%
\[
\mathcal{H}(x,x^{\ast})\geq\left\langle x,x^{\ast}\right\rangle .
\]
Moreover, the following assertions are equivalent:$\smallskip$\newline$%
\begin{array}
[c]{l}%
\left(  a\right)  \quad(x,x^{\ast})\in A;\smallskip\\
\left(  b\right)  \quad\mathcal{H}(x,x^{\ast})=\left\langle x,x^{\ast
}\right\rangle ;\smallskip\\
\left(  c\right)  \quad\mathcal{H}^{\ast}(x^{\ast},x)=\left\langle x,x^{\ast
}\right\rangle ;\smallskip\\
\left(  d\right)  \quad\exists~\left(  u,u^{\ast}\right)  \in Dom\left(
\partial\mathcal{H}\right)  \text{ such that }\left(  x^{\ast},x\right)
\in\partial\mathcal{H}\left(  u,u^{\ast}\right)  \ \text{and}\ \left\langle
u-x,u^{\ast}-x^{\ast}\right\rangle =0;\smallskip\\
\left(  e\right)  \quad\left(  x^{\ast},x\right)  \in\partial\mathcal{H}%
\left(  x,x^{\ast}\right)  .
\end{array}
$
\end{theorem}

\begin{proof}
It is not difficult to show that $\left(  b\right)  \Leftrightarrow\left(
a\right)  \Rightarrow\left(  e\right)  \Rightarrow\left(  d\right)
\Rightarrow\left(  a\right)  $. Moreover, using the Fenchel equality:%
\[
\left(  x^{\ast},x\right)  \in\partial\mathcal{H}\left(  x,x^{\ast}\right)
\;\Rightarrow\;\mathcal{H}(x,x^{\ast})+\mathcal{H}^{\ast}(x^{\ast
},x)=\left\langle (x,x^{\ast}),(x^{\ast},x)\right\rangle ,
\]
we obtain that $\left(  e\right)  \&\left(  b\right)  $\ $\Rightarrow\ \left(
c\right)  .$ The point $\left(  c\right)  $ yields $\left(  a\right)  $ by
using the equivalent form of the definition of $\mathcal{H}^{\ast}:$%
\[
\mathcal{H}^{\ast}(x^{\ast},x)=\left\langle x,x^{\ast}\right\rangle
-\inf_{\left(  u,u^{\ast}\right)  \in\mathbb{X}\times\mathbb{X}^{\ast}%
}\left\{  \left\langle x-u,x^{\ast}-u^{\ast}\right\rangle +\mathcal{H}\left(
u,u^{\ast}\right)  -\left\langle u,u^{\ast}\right\rangle \right\}  .
\]

\hfill
\end{proof}

\begin{remark}
The function $\mathcal{H}_{A}$ is minimal in the family of convex functions
$f:\mathbb{X}\times\mathbb{X}^{\ast}\rightarrow\mathbb{]-}\infty,+\infty]$
with the properties: $f(x,x^{\ast})\geq\left\langle x,x^{\ast}\right\rangle
\;$for all $(x,x^{\ast})\in\mathbb{X}\times\mathbb{X}^{\ast}$ and
$f(x,x^{\ast})=\left\langle x,x^{\ast}\right\rangle \;$for all $(x,x^{\ast
})\in A.$
\end{remark}

Using the above tools, in the paper \cite{SZ}, Simons and Z\u{a}linescu give a
nice proof of the famous Rockafellar's characterization of a maximal monotone
operator.\medskip

Let $\mathbb{H}$ be a real separable Hilbert space and $A:\mathbb{H}%
\rightrightarrows\mathbb{H}$ be a maximal monotone operator. Denote for
$\varepsilon>0,$ $J_{\varepsilon},A_{\varepsilon}:\mathbb{H}\rightarrow
\mathbb{H},$ the ($1$-, resp. $1/\varepsilon~$-) Lipschitz continuous
functions $J_{\varepsilon}\left(  x\right)  =\left(  I+\varepsilon A\right)
^{-1}\left(  x\right)  $ and%
\[
A_{\varepsilon}\left(  x\right)  =\frac{x-J_{\varepsilon}\left(  x\right)
}{\varepsilon}\in A\left(  J_{\varepsilon}\left(  x\right)  \right)  .
\]

Let%
\[
BV_{0}\left(  \left[  0,T\right]  ;\mathbb{H}\right)  =\left\{  k:\left[
0,T\right]  \rightarrow\mathbb{H}:\left\updownarrow k\right\updownarrow
_{T}<\infty,\;k\left(  0\right)  =0\right\}  ,
\]
where $\left\updownarrow k\right\updownarrow _{T}:=\left\Vert k\right\Vert
_{BV\left(  \left[  0,T\right]  ;\mathbb{H}\right)  }.$ If we consider on
$C\left(  \left[  0,T\right]  ;\mathbb{H}\right)  $ the usual norm%
\[
\left\Vert y\right\Vert _{C\left(  \left[  0,T\right]  ;\mathbb{H}\right)
}=\left\Vert y\right\Vert _{T}=\sup\left\{  \left\vert y\left(  s\right)
\right\vert :0\leq s\leq T\right\}  ,
\]
then $\left(  C\left(  \left[  0,T\right]  ;\mathbb{H}\right)  \right)
^{\ast}=BV_{0}\left(  \left[  0,T\right]  ;\mathbb{H}\right)  .$ We denote the
duality between these spaces by%
\[
\left\langle \!\left\langle z,g\right\rangle \!\right\rangle :=\int_{0}%
^{T}\left\langle z\left(  t\right)  ,dg\left(  t\right)  \right\rangle .
\]

Denote by $\mathcal{A}$ the realization on $C\left(  \left[  0,T\right]
;\mathbb{H}\right)  $ of the maximal monotone operator\newline$A:\mathbb{H}%
\rightrightarrows\mathbb{H},$ that is the operator $\mathcal{A}:C\left(
\left[  0,T\right]  ;\mathbb{H}\right)  \rightrightarrows BV_{0}\left(
\left[  0,T\right]  ;\mathbb{H}\right)  $ defined as follows: $\left(
x,k\right)  \in\mathcal{A}$ if $x\in C\left(  \left[  0,T\right]
;\mathbb{R}^{d}\right)  $, $k\in BV_{0}\left(  \left[  0,T\right]
;\mathbb{H}\right)  $ and one of the following equivalent conditions are satisfied:

\begin{itemize}
\item[$\left(  d_{1}\right)  $] for all $0\leq s\leq t\leq T,$ $%
{\displaystyle\int_{s}^{t}}
\left\langle x\left(  r\right)  -z,dk\left(  r\right)  -z^{\ast}%
dr\right\rangle \geq0,\ \forall\left(  z,z^{\ast}\right)  \in A;$

\item[$\left(  d_{2}\right)  $] for all $0\leq s\leq t\leq T$ and for all
$u,u^{\ast}\in C([0,T];\mathbb{H})\ $such that $\left(  u(r),u^{\ast
}(r)\right)  \in A,\ \forall r\in\lbrack s,t],$%
\[%
{\displaystyle\int\nolimits_{s}^{t}}
\left\langle x(r)-u(r),dk(r)-u^{\ast}(r)dr\right\rangle \geq0;
\]

\item[$\left(  d_{3}\right)  $] for all $u,u^{\ast}\in C([0,T];\mathbb{H}%
)\ $such that $\left(  u(r),u^{\ast}(r)\right)  \in A,$ $\forall r\in\left[
0,T\right]  ,$%
\[%
{\displaystyle\int_{0}^{T}}
\left\langle x(r)-u(r),dk(r)-u^{\ast}(r)dr\right\rangle \geq0.
\]

\end{itemize}

\noindent$\mathcal{A}$ is a maximal monotone operator since, setting%
\[
u\left(  r\right)  =J_{\varepsilon}\left(  \frac{x\left(  r\right)  +y\left(
r\right)  }{2}\right)  =\frac{x\left(  r\right)  +y\left(  r\right)  }%
{2}-\varepsilon A_{\varepsilon}\left(  \frac{x\left(  r\right)  +y\left(
r\right)  }{2}\right)  \text{\ };\ u^{\ast}\left(  r\right)  =A_{\varepsilon
}\left(  \frac{x\left(  r\right)  +y\left(  r\right)  }{2}\right)
\]
in $\left(  d_{2}\right)  $ written for $\left(  x,k\right)  \in\mathcal{A}$
and respectively for $(y,\ell)\in\mathcal{A}$ and taking then $\varepsilon
\rightarrow0,$ we infer (since $\varepsilon A_{\varepsilon}\rightarrow0$ as
$\varepsilon\rightarrow0$) that%
\begin{equation}
\int_{s}^{t}\left\langle x\left(  r\right)  -y\left(  r\right)  ,dk\left(
r\right)  -d\ell\left(  r\right)  \right\rangle \geq0,\ \ \forall0\leq s\leq
t\leq T. \label{pozitivitate max mon}%
\end{equation}
The maximality clearly follows from the definition of $\mathcal{A}$.

For the realization of the operator $A$ on $L^{r}\left(  0,T;\mathbb{H}%
\right)  ,$ $r\geq1,$ we use the same notation $\mathcal{A}$ without risk of
confusion since every time we mention the space of realization. In this case,
the operator $\mathcal{A}:L^{r}\left(  0,T;\mathbb{H}\right)
\rightrightarrows L^{q}\left(  0,T;\mathbb{H}\right)  ,$ $\dfrac{1}{r}%
+\dfrac{1}{q}=1$ is defined by $\left(  x,g\right)  \in\mathcal{A}$ if

\begin{itemize}
\item $%
{\displaystyle\int_{s}^{t}}
\left\langle x\left(  r\right)  -z,g\left(  r\right)  -z^{\ast}\right\rangle
dr\geq0,\ $for all $0\leq s\leq t\leq T$ and for all $\left(  z,z^{\ast
}\right)  \in A,$
\end{itemize}

or (clearly), equivalently

\begin{itemize}
\item $g\left(  t\right)  \in A\left(  x\left(  t\right)  \right)  ,\;a.e.$
$t\in\left[  0,T\right]  .$
\end{itemize}

\noindent Arguing similar to the previous situation, we obtain that
$\mathcal{A}$ is a maximal monotone operator.

\subsection{Stochastic framework\label{sfram}}

Let $({\Omega},\mathcal{F},\mathbb{P},\{\mathcal{F}_{t}\}_{t\geq0})$ be a
stochastic basis i.e. $({\Omega},\mathcal{F},\mathbb{P})$ is a complete
probability space and $\{\mathcal{F}_{t}\}_{t\geq0}$ is a filtration
satisfying the usual assumptions of right continuity and completeness:%
\[
\mathcal{N}_{\mathbb{P}}\subset\mathcal{F}_{s}\subset\mathcal{F}_{t}=%
{\textstyle\bigcap\nolimits_{\varepsilon>0}}
\mathcal{F}_{t+\varepsilon}~,
\]
for all $0\leq s\leq t$, where $\mathcal{N}_{\mathbb{P}}$ is the set of all
$\mathbb{P}$-null sets.

Let $\left(  \mathbb{H},\left\vert ~\cdot~\right\vert _{\mathbb{H}}\right)  $
be a real separable Hilbert space; if $F$ is a closed subset of $\mathbb{H}$,
denote by $\mathcal{B}_{F}$ the $\sigma$-algebra generated by the closed
subsets of $F.$

Denote by $S_{\mathbb{H}}^{p}\left[  0,T\right]  ,$ $p\geq0,$ the space of
progressively measurable continuous stochastic processes $X:\Omega
\times\left[  0,T\right]  \rightarrow\mathbb{H}$ (i.e. $t\longmapsto X\left(
\omega,t\right)  $ is continuous a.s. $\omega\in\Omega$, and $\left(
\omega,s\right)  \longmapsto X\left(  \omega,s\right)  :\Omega\times\left[
0,T\right]  \rightarrow\mathbb{H}$ is $\left(  \mathcal{F}_{t}\otimes
\mathcal{B}_{\left[  0,t\right]  },\mathcal{B}_{\mathbb{H}}\right)  $
measurable for all $t\in\left[  0,T\right]  $), such that%
\[
\left\Vert X\right\Vert _{S_{\mathbb{H}}^{p}\left[  0,T\right]  }=\left\{
\begin{array}
[c]{ll}%
\left(  \mathbb{E~}\left\Vert X\right\Vert _{T}^{p}\right)  ^{\frac{1}%
{p}\wedge1}<{\infty}, & \;\text{if }p>0,\bigskip\\
\mathbb{E}\left[  1\wedge\left\Vert X\right\Vert _{T}\right]  , & \;\text{if
}p=0,
\end{array}
\right.
\]
where%
\[
\left\Vert X\right\Vert _{T}:=\sup_{t\in\left[  0,T\right]  }\left\vert
X_{t}\right\vert .
\]
The space $(S_{\mathbb{H}}^{p}\left[  0,T\right]  ,\left\Vert \cdot\right\Vert
_{S_{\mathbb{H}}^{p}\left[  0,T\right]  }),\ p\!\geq1,$ is a Banach space and
$S_{\mathbb{H}}^{p}\left[  0,T\right]  $, $0\leq p<1,$ is a complete metric
space with the metric $\rho(Z_{1},Z_{2})\!=\left\Vert Z_{1}\!-\!Z_{2}%
\right\Vert _{S_{d}^{p}\left[  0,T\right]  }$ (when $p=0$ the metric
convergence coincides with the probability convergence).

If $\mathbb{H=R}^{d}$ we will denote $S_{\mathbb{H}}^{p}\left[  0,T\right]  $
by $S_{d}^{p}\left[  0,T\right]  $.

Let $\left(  \mathbb{H}_{0},\left\vert \cdot\right\vert _{\mathbb{H}_{0}%
}\right)  $ be a real separable Hilbert space and%
\[
B=\{B_{t}(\varphi):\left(  t,\varphi\right)  \in\lbrack0,T]\times
\mathbb{H}_{0}\}\subset L^{0}\left(  \Omega,\mathcal{F},\mathbb{P}\right)
\]
a Gaussian family of real-valued random variables with zero mean and
covariance function%
\[
\mathbb{E}\left[  B_{t}(\varphi)B_{s}(\psi)\right]  =\left(  t\wedge s\right)
\times\left\langle \varphi,\psi\right\rangle _{\mathbb{H}_{0}}\,,\quad
\forall\varphi,\psi\in\mathbb{H}_{0},\quad\forall s,t\in\lbrack0,T],
\]
where $t\wedge s=\min\left\{  t,s\right\}  .$ We call $\left(  B{,}%
\{\mathcal{F}_{t}\}\right)  $ a $\mathbb{H}_{0}$-\textit{Wiener process }if,
for all $t\in\lbrack0,T]$, we have%
\[%
\begin{array}
[c]{l}%
\ \left(  i\right)  \quad\mathcal{F}_{t}^{B}=\sigma\{B_{s}(\varphi
);\;s\in\lbrack0,t],\varphi\in\mathbb{H}_{0}\}\vee\mathcal{N}_{\mathbb{P}%
}\;\subset\mathcal{F}_{t}\text{ and}\smallskip\\
\left(  ii\right)  \quad B_{t+h}(\varphi)-B_{t}(\varphi)\text{ is independent
of }\mathcal{F}_{t}\,\text{, for all }h>0,\text{ }\varphi\in\mathbb{H}_{0}.
\end{array}
\]
Note that, given any orthonormal basis $\{e_{i};\;i\in I\subseteq
\mathbb{N}^{\ast}\}$ of $\mathbb{H}_{0},$ the sequence $\beta^{i}=\{\beta
_{t}^{i}=B_{t}(e_{i});\;t\in\lbrack0,T]\},\;i\in I$, defines a family of
independent real-valued standard Wiener processes (Brownian motions).
Moreover, if $\mathbb{H}_{0}$ is of finite dimension, we have%
\[
B_{t}=\sum\limits_{i\geq1}\beta_{t}^{i}e_{i}\,,\quad t\in\lbrack0,T]\text{.}%
\]
In the general case this series does not converge in $\mathbb{H}_{0}$, but
rather in a larger space $\mathbb{\tilde{H}}_{0},\newline\mathbb{H}_{0}%
\subset\mathbb{\tilde{H}}_{0}$ which is such that the injection of
$\mathbb{H}_{0}$ into $\mathbb{\tilde{H}}_{0}$ is Hilbert-Schmidt.
Moreover,\newline$B\in\mathcal{M}^{2}(0,T;\mathbb{\tilde{H}}_{0})$.

By $\mathcal{M}^{p}(0,T;\mathbb{H}),\ p\geq1,$ we denote the space of
$\mathbb{H}$-valued continuous, $p$-integrable martingales $M$, that is, the
space of all continuous stochastic processes $M:\Omega\times\left[
0,T\right]  \rightarrow\mathbb{H}$ satisfying, $\mathbb{P}$-$a.s$,%
\[%
\begin{array}
[c]{ll}%
(m_{1})\quad & M_{0}=0,\smallskip\\
(m_{2})\quad & \mathbb{E}\left\vert M_{t}\right\vert ^{p}<\infty,\ \forall
t\in\left[  0,T\right]  ,\smallskip\\
(m_{3})\quad & \mathbb{E}\left[  M_{t}|\mathcal{F}_{s}\right]  =M_{s},\text{
for all }0\leq s\leq t\leq T.
\end{array}
\]
$\mathcal{M}^{p}(0,T;\mathbb{H})$ is a Banach space with respect to the norm
$\left\Vert X\right\Vert _{\mathcal{M}^{p}}=\left(  \mathbb{E~}\left\vert
X_{T}\right\vert ^{p}\right)  ^{1/p}$; in the case $p>1$, $\mathcal{M}%
^{p}(0,T;\mathbb{H})$ is a closed linear subspace of $S_{\mathbb{H}}%
^{p}\left[  0,T\right]  $.

In order to define the stochastic integral with respect to the $\mathbb{H}%
_{0}$-Wiener process $B$, we introduce a class of processes with values in the
separable Hilbert space $\mathcal{L}^{2}(\mathbb{H}_{0};\mathbb{H})$ of
Hilbert--Schmidt operators from $\mathbb{H}_{0}$ into $\mathbb{H}$,
\textit{i.e}. the space of linear operators $F:\mathbb{H}_{0}\rightarrow
\mathbb{H}$ satisfying%
\[
\left\Vert F\right\Vert _{HS}^{2}=\sum_{i=1}^{\infty}\left\vert Fe_{i}%
\right\vert _{\mathbb{H}}^{2}=\mathbf{Tr}F^{\ast}F=\mathbf{Tr}FF^{\ast}%
<\infty.
\]
Denote $\Lambda_{\mathbb{H}\times\mathbb{H}_{0}}^{p}\left(  0,T\right)
,\ p\in\lbrack0,{\infty}[,$ the space of progressively measurable
processes\newline$Z:{\Omega}\times]0,T[\rightarrow\mathcal{L}^{2}%
(\mathbb{H}_{0};\mathbb{H})$ such that:%
\[
\left\Vert Z\right\Vert _{\Lambda^{p}}=\left\{
\begin{array}
[c]{ll}%
\left[  \mathbb{E~}\left(
{\displaystyle\int_{0}^{T}}
\Vert Z_{s}\Vert_{HS}^{2}ds\right)  ^{\frac{p}{2}}\right]  ^{\frac{1}{p}%
\wedge1}, & \;\text{if }p>0,\bigskip\\
\mathbb{E}\left[  1\wedge\left(
{\displaystyle\int_{0}^{T}}
\Vert Z_{s}\Vert_{HS}^{2}ds\right)  ^{\frac{1}{2}}\right]  , & \;\text{if
}p=0.
\end{array}
\right.
\]
The space $(\Lambda_{\mathbb{H}\times\mathbb{H}_{0}}^{p}\left(  0,T\right)
,\left\Vert \cdot\right\Vert _{\Lambda^{p}}),\ p\geq1,$ is a Banach space and
$\Lambda_{\mathbb{H}\times\mathbb{H}_{0}}^{p}\left(  0,T\right)  $, $0\leq
p<1,$ is a complete metric space with the metric $\rho(Z_{1},Z_{2})=\left\Vert
Z_{1}-Z_{2}\right\Vert _{\Lambda^{p}}$.

Consider $\{e_{i};\ i\in I\subset\mathbb{N}^{\ast}\}$ an orthonormal basis of
$\mathbb{H}_{0}.$ Let $Z\in\Lambda_{\mathbb{H}\times\mathbb{H}_{0}}^{p}\left(
0,T\right)  $, with $p\geq0$. The stochastic integral $I$ is defined by
$Z\overset{I}{\mapsto}I_{\cdot}\left(  Z\right)  $, where%
\[
I_{t}\left(  Z\right)  :=\int_{0}^{t}Z_{s}dB_{s}=\sum\limits_{i\in I}\int
_{0}^{t}Z_{s}(e_{i})dB_{s}(e_{i}),\text{ }t\in\left[  0,T\right]  .
\]
Note that it doesn't depend on the choice of the orthonormal basis of
$\mathbb{H}_{0}$. The application%
\[
I:\Lambda_{\mathbb{H}\times\mathbb{H}_{0}}^{p}\left(  0,T\right)  \rightarrow
S_{\mathbb{H}}^{p}\left[  0,T\right]
\]
is a linear continuous operator and it has the following properties:%

\[%
\begin{array}
[c]{ll}%
(a)\quad & \mathbb{E}I_{t}\left(  Z\right)  =0,\quad\text{if }p\geq
1,\vspace*{1mm}\\
(b)\quad & \mathbb{E}|I_{T}\left(  Z\right)  |^{2}=\left\Vert Z\right\Vert
_{\Lambda^{2}}^{2}~,\quad\text{if }p\geq2,\vspace*{1mm}\\
(c)\quad & \dfrac{1}{c_{p}}\left\Vert Z\right\Vert _{\Lambda^{p}}^{p}%
\leq\mathbb{E}\sup\limits_{t\in\lbrack0,T]}|I_{t}\left(  Z\right)  |^{p}\leq
c_{p}\left\Vert Z\right\Vert _{\Lambda^{p}}^{p}~,\text{ if }p>0,\smallskip\\
& \quad\quad\text{(Burkholder-Davis-Gundy inequality)}\vspace*{1mm}%
\smallskip\\
(d)\quad & I(Z)\in\mathcal{M}^{p}(\Omega\times\lbrack0,T];\mathbb{H}),\text{
}p\geq1.
\end{array}
\]
The definition and the properties of the stochastic integral can be found in
Pardoux \& R\u{a}\c{s}canu \cite{PRb} or Da Prato \& Zabczyk \cite{DZ}.

If $\mathbb{H}_{0}=\mathbb{R}^{k}$ and $\mathbb{H=R}^{d}$ then $\{B_{t}%
,t\geq0\}$ is a $k$-dimensional Wiener process (Brownian motion);
$\mathcal{L}^{2}(\mathbb{H}_{0};\mathbb{H})$ is the space of real matrices
$F=\left(  f_{ij}\right)  _{d\times k}$ and $\left\vert F\right\vert
^{2}:=\left\Vert F\right\Vert _{HS}^{2}=%
{\displaystyle\sum\limits_{i,j}}
f_{i,j}^{2}$. In this situation, the space $\Lambda_{\mathbb{H}\times
\mathbb{H}_{0}}^{p}\left(  0,T\right)  $ will be denoted by $\Lambda_{d\times
k}^{p}\left(  0,T\right)  .$

\section{Fitzpatrick function approach}

\subsection{A Generalized Skorohod problem}

Throughout this section $\mathbb{H}$ is a real separable Hilbert space with
the norm $\left\vert \cdot\right\vert $ and the scalar product $\left\langle
\cdot,\cdot\right\rangle $.

We study the multivalued monotone differential equation%
\begin{equation}
\left\{
\begin{array}
[c]{l}%
dx\left(  t\right)  +Ax\left(  t\right)  (dt)\ni dm\left(  t\right)
,\medskip\\
x\left(  0\right)  =x_{0},\quad t\geq0,
\end{array}
\right.  \quad\quad\left(  GSP\right)  \label{3}%
\end{equation}
where we assume%
\[
\left(  H_{GSP}\right)  :\quad\left\{
\begin{array}
[c]{rl}%
\left(  i\right)  & A:\mathbb{H}\rightrightarrows\mathbb{H}\text{ is a maximal
monotone operator,}\medskip\\
\left(  ii\right)  & x_{0}\in\overline{Dom(A)},\medskip\\
\left(  iii\right)  & m:[0,\infty)\longrightarrow\mathbb{H}\text{ is
continuous and }m\left(  0\right)  =0.
\end{array}
\right.
\]

\begin{definition}
A continuous function $x:\left[  0,T\right]  \rightarrow\mathbb{H}$ is a
solution of Eq.(\ref{3}) if $x\left(  t\right)  \in\overline{Dom(A)}$ for all
$0\leq t\leq T,$ ($T$ arbitrarily fixed) and there exists $k\in C\left(
\left[  0,T\right]  ;\mathbb{H}\right)  \bigcap BV_{0}\left(  \left[
0,T\right]  ;\mathbb{H}\right)  $ such that%
\[
x\left(  t\right)  +k\left(  t\right)  =x_{0}+m\left(  t\right)
,\ \forall0\leq t\leq T
\]
and%
\begin{equation}%
{\displaystyle\int_{s}^{t}}
\left\langle x\left(  r\right)  -z,dk\left(  r\right)  -z^{\ast}%
dr\right\rangle \geq0,\ \forall\left(  z,z^{\ast}\right)  \in A,\ \forall0\leq
s\leq t\leq T. \label{integr def sol GSP}%
\end{equation}
(Without confusion, the uniqueness of $k$ will permit us to call the pair
$\left(  x,k\right)  $ solution of the generalized Skorohod problem $(GSP)$
and we write $\left(  x,k\right)  =\mathcal{GSP}\left(  A;x_{0},m\right)  $.)
\end{definition}

In virtue of this definition, the (classical) Skorohod problem (for more
details, one can consult C\'{e}pa \cite{C1} or \cite{C2}) is obtained for
$A=\partial I_{E}:\mathbb{R}^{d}\rightrightarrows\mathbb{R}^{d},$ where $E$ is
a closed convex subset of $\mathbb{R}^{d},$%
\[
I_{E}\left(  x\right)  =\left\{
\begin{array}
[c]{rl}%
0, & \text{if }x\in E,\\
+\infty, & \text{if }x\in\mathbb{R}^{d}\setminus E
\end{array}
\right.
\]
and%
\[
\partial I_{E}\left(  x\right)  =\left\{
\begin{array}
[c]{ll}%
0, & \text{if }x\in int(E),\\
\left\{  \nu\in\mathbb{R}^{d}:\left\langle \nu,y-x\right\rangle \leq0,\text{
for all }y\in E\right\}  , & \text{if }x\in Bd\left(  E\right)  ,\\
\emptyset, & \text{if }x\notin E.
\end{array}
\right.
\]
The definition of the solution can be given in a equivalent form as follows.

\begin{definition}
A continuous function $x:\left[  0,T\right]  \rightarrow\mathbb{R}^{d}%
\times\mathbb{R}^{d}$ is a solution of Skorohod problem in $E$ if $x\left(
t\right)  \in E$ for all $0\leq t\leq T$ and there exists $k\in C\left(
\left[  0,T\right]  ;\mathbb{R}^{d}\right)  \bigcap BV_{0}\left(  \left[
0,T\right]  ;\mathbb{R}^{d}\right)  $ such that%
\[
\left\{
\begin{array}
[c]{ll}%
\left(  a\right)  \quad & \left\updownarrow k\right\updownarrow _{t}=%
{\displaystyle\int_{0}^{t}}
\mathbf{1}_{x\left(  s\right)  \in Bd\left(  E\right)  }d\left\updownarrow
k\right\updownarrow _{s},\medskip\\
\left(  b\right)  \quad &
\begin{array}
[c]{l}%
k\left(  t\right)  =%
{\displaystyle\int_{0}^{t}}
n_{x\left(  s\right)  }d\left\updownarrow k\right\updownarrow _{s}%
,\;\text{where }n_{x\left(  s\right)  }\in N_{E}\left(  x\left(  s\right)
\right) \\
\quad\quad\quad\quad\text{and }\left\vert n_{x\left(  s\right)  }\right\vert
=1,\;d\updownarrow k\updownarrow_{s}-a.e.
\end{array}
\end{array}
\right.
\]
and%
\[
x\left(  t\right)  +k\left(  t\right)  =x_{0}+m\left(  t\right)  ,\ \forall
t\in\left[  0,T\right]  .
\]
($N_{E}\left(  x\right)  $ denotes the outward normal cone to $E$ at $x\in E$.)
\end{definition}

\medskip

Let $\mathcal{A}:C\left(  \left[  0,T\right]  ;\mathbb{H}\right)
\rightrightarrows BV_{0}\left(  \left[  0,T\right]  ;\mathbb{H}\right)  $ be
the realization of the maximal monotone operator $A:\mathbb{H}%
\rightrightarrows\mathbb{H}$ and%
\[
\mathbb{X}=\left\{  \mu\in C\left(  \left[  0,T\right]  ;\mathbb{H}\right)
:\mu\left(  0\right)  =0\right\}
\]
the linear closed subspace of $C\left(  \left[  0,T\right]  ;\mathbb{H}%
\right)  .$ For each $R>0$, we define%
\[
\mathbb{Y}_{R}=\left\{  k\in C\left(  \left[  0,T\right]  ;\mathbb{H}\right)
:k\left(  0\right)  =0,\;\left\updownarrow k\right\updownarrow _{T}\leq
R\right\}  ;
\]
$\mathbb{Y}_{R}$ is a closed subset of $C\left(  \left[  0,T\right]
;\mathbb{H}\right)  $ and, consequently, it is a metric space with respect to
the metric from $C\left(  \left[  0,T\right]  ;\mathbb{H}\right)  $. Remark
that, by Helly-Foia\c{s} theorem (see Barbu \& Precupanu \cite{BP}, Theorem
3.5 \& Remark 3.2), it is also a bounded $w^{\ast}$-closed subset of
$BV_{0}\left(  \left[  0,T\right]  ;\mathbb{H}\right)  $.

Let $\alpha:\mathbb{R}_{+}\rightarrow\mathbb{R}_{+}$ a continuous function
such that $\alpha\left(  0\right)  =0.$ Denote%
\[
C_{\alpha}=\left\{  x\in\mathbb{X}:\mathbf{m}_{x}\left(  \varepsilon\right)
\leq\alpha\left(  \varepsilon\right)  \;\text{for all }\varepsilon
\geq0\right\}  .
\]
Here the function $\mathbf{m}_{x}:\mathbb{R}_{+}\rightarrow\mathbb{R}_{+}$
represents the \textit{modulus of continuity} of the continuous function
$x:\left[  0,T\right]  \rightarrow\mathbb{H}$ and it is defined by%
\[
\mathbf{m}_{x}\left(  \delta\right)  =\mathbf{m}_{x,T}\left(  \delta\right)
=\sup\left\{  \left\vert x\left(  t\right)  -x\left(  s\right)  \right\vert
:\left\vert t-s\right\vert \leq\delta,\;t,s\in\left[  0,T\right]  \right\}  .
\]
Clearly, $C_{\alpha}$ is a bounded closed convex subset of $\mathbb{X}%
.\smallskip$

Consider, for each $\left(  u,u^{\ast}\right)  \in\mathcal{A}$ and $\nu
\in\mathbb{X}$, the function $J_{\left(  u,u^{\ast},\nu\right)  }%
:\mathbb{H}\times\mathbb{X}\times\mathbb{Y}_{R}\times\mathbb{X}\rightarrow
\mathbb{R}$ given by%
\begin{align*}
J_{\left(  u,u^{\ast},\nu\right)  }\left(  a,x,k,\mu\right)   &  =\left\vert
a-x_{0}\right\vert ^{2}+%
{\displaystyle\int_{0}^{T}}
\left[  \left\langle u\left(  t\right)  ,dk\left(  t\right)  \right\rangle
+\left\langle x\left(  t\right)  ,du^{\ast}\left(  t\right)  \right\rangle
-\left\langle u\left(  t\right)  ,du^{\ast}\left(  t\right)  \right\rangle
\right] \\
&  -\int_{0}^{T}\left\langle x\left(  t\right)  ,dk\left(  t\right)
\right\rangle +2R\left\Vert \mu-m\right\Vert _{T}+%
{\displaystyle\int_{0}^{T}}
\left\langle \mu\left(  t\right)  -\nu\left(  t\right)  ,dk\left(  t\right)
\right\rangle -R\left\Vert \nu-m\right\Vert _{T}%
\end{align*}
and $\hat{J}:\mathbb{H}\times\mathbb{X}\times\mathbb{Y}_{R}\times
\mathbb{X}\rightarrow]-\infty,+\infty],$ defined by%
\begin{equation}%
\begin{array}
[c]{l}%
\hat{J}\left(  a,x,k,\mu\right)  =\underset{\left(  u,u^{\ast}\right)
\in\mathcal{A},\ \nu\in C_{\alpha}}{\sup}J_{\left(  u,u^{\ast},\nu\right)
}\left(  a,x,k,\mu\right)  \medskip\\
\quad=\left\vert a-x_{0}\right\vert ^{2}+\mathcal{H}\left(  x,k\right)
-\left\langle \!\left\langle x,k\right\rangle \!\right\rangle +2R\left\Vert
\mu-m\right\Vert _{T}+\sup\limits_{\nu\in C_{\alpha}}\left\{  \left\langle
\!\left\langle \mu-\nu,k\right\rangle \!\right\rangle -R\left\Vert
\nu-m\right\Vert _{T}\right\}  ,
\end{array}
\label{fcJ}%
\end{equation}
where $\mathcal{H}:C\left(  \left[  0,T\right]  ;\mathbb{H}\right)  \times
BV_{0}\left(  \left[  0,T\right]  ;\mathbb{H}\right)  \rightarrow
]-\infty,+\infty]$ is the Fitzpatrick function associated to the maximal
monotone operator $\mathcal{A}$.

\begin{remark}
$\hat{J}:\mathbb{H}\times\mathbb{X}\times\mathbb{Y}_{R}\times\mathbb{X}%
\rightarrow]-\infty,+\infty]$ is a lower semicontinuous function as the
supremum of the continuous functions $J_{\left(  u,u^{\ast},\nu\right)  }$.
\end{remark}

\noindent Remark also that, for $\mu\in C_{\alpha}$,%
\[
2R\left\Vert \mu-m\right\Vert _{T}+\sup_{\nu\in C_{\alpha}}\left\{
\left\langle \!\left\langle \mu-\nu,k\right\rangle \!\right\rangle
-R\left\Vert \nu-m\right\Vert _{T}\right\}  \geq R\left\Vert \mu-m\right\Vert
_{T}\geq0.
\]

\begin{proposition}
\label{p1}Let $R>0$ and $\alpha:\mathbb{R}_{+}\rightarrow\mathbb{R}_{+}$ a
continuous function such that $\alpha\left(  0\right)  =0.$ The function
$\hat{J}$ has the following properties

\begin{itemize}
\item[$\left(  a\right)  $] $\;\hat{J}\left(  a,x,k,\mu\right)  \geq0$, for
all $\left(  a,x,k,\mu\right)  \in\mathbb{H}\times\mathbb{X}\times
\mathbb{Y}_{R}\times C_{\alpha}.$

\item[$\left(  b\right)  $] $\;$Let $(\hat{a},\hat{x},\hat{k},\hat{\mu}%
)\in\mathbb{H}\times\mathbb{X}\times\mathbb{Y}_{R}\times C_{\alpha}.$ Then
$\hat{J}(\hat{a},\hat{x},\hat{k},\hat{\mu})=0$ iff $\hat{a}=x_{0},$ $\hat{\mu
}=m$ and $\hat{k}\in\mathcal{A}\left(  \hat{x}\right)  .$

\item[$\left(  c\right)  $] $\;$The restriction of $\hat{J}$ to the closed
convex set%
\[
\mathbb{K}=\left\{  \left(  a,x,k,\mu\right)  \in\mathbb{H}\times
\mathbb{X}\times\mathbb{Y}_{R}\times C_{\alpha}:x+k=a+\mu\right\}
\]
is a convex lower semicontinuous function; for $(\hat{a},\hat{x},\hat{k}%
,\hat{\mu})\in\mathbb{K},$ we have%
\[
\hat{J}(\hat{a},\hat{x},\hat{k},\hat{\mu})=0\quad\text{iff}\quad\hat{a}%
=x_{0},\ \hat{\mu}=m\text{ and }(\hat{x},\hat{k})=\mathcal{GSP}\left(
A;x_{0},m\right)  .
\]

\end{itemize}
\end{proposition}

\begin{proof}
The points $\left(  a\right)  $ and $\left(  b\right)  $ clearly are
consequences of the properties of the Fitzpatrick function $\mathcal{H}.$ Let
us prove $\left(  c\right)  .$ We have $\left(  a,x,k,\mu\right)
\in\mathbb{K}$ and%
\begin{align*}
\hat{J}\left(  a,x,k,\mu\right)   &  =\left\vert a-x_{0}\right\vert
^{2}+\mathcal{H}\left(  x,k\right)  -\left\langle \!\left\langle
x,k\right\rangle \!\right\rangle +2R\left\Vert \mu-m\right\Vert _{T}\medskip\\
&  +\sup_{\nu\in C_{\alpha}}\left\{  \left\langle \!\left\langle \mu
-\nu,k\right\rangle \!\right\rangle -R\left\Vert \nu-m\right\Vert _{T}\right\}
\\
&  =\left\vert a-x_{0}\right\vert ^{2}+\mathcal{H}\left(  x,k\right)
+\frac{1}{2}\left\vert x\left(  T\right)  -\mu\left(  T\right)  \right\vert
^{2}-\frac{1}{2}\left\vert a\right\vert ^{2}-%
{\displaystyle\int_{0}^{T}}
\left\langle \mu\left(  s\right)  ,dk\left(  s\right)  \right\rangle
\medskip\\
&  +2R\left\Vert \mu-m\right\Vert _{T}+\sup_{\nu\in C_{\alpha}}\left\{
\left\langle \!\left\langle \mu-\nu,k\right\rangle \!\right\rangle
-R\left\Vert \nu-m\right\Vert _{T}\right\} \\
&  =\left\vert x_{0}\right\vert ^{2}-2\left\langle a,x_{0}\right\rangle
+\frac{1}{2}\left\vert a\right\vert ^{2}+\mathcal{H}\left(  x,k\right)
+\frac{1}{2}\left\vert x\left(  T\right)  -\mu\left(  T\right)  \right\vert
^{2}\medskip\\
&  +2R\left\Vert \mu-m\right\Vert _{T}+\sup_{\nu\in C_{\alpha}}\left\{
\left\langle \!\left\langle -\nu,k\right\rangle \!\right\rangle -R\left\Vert
\nu-m\right\Vert _{T}\right\}
\end{align*}
and the convexity of $\hat{J}$ follows.\hfill
\end{proof}

\bigskip

In the sequel we prove the existence and uniqueness of the solution of the
multivalued monotone differential equation (\ref{3}). Our proof is strongly
connected with the one from R\u{a}\c{s}canu \cite{R}. First highlight some
properties of a solution $\left(  x,k\right)  =\mathcal{GSP}\left(
A;x_{0},m\right)  .$

Consider $\mathcal{M}$ a bounded and equicontinuous subset of $C\left(
\left[  0,T\right]  ;\mathbb{H}\right)  $ and we denote%
\[
\left\Vert \mathcal{M}\right\Vert _{T}=\sup\left\{  \left\Vert y\right\Vert
_{T}:y\in\mathcal{M}\right\}  \quad\text{and}\quad\mathbf{m}_{\mathcal{M}%
,T}\left(  \delta\right)  =\sup\left\{  \mathbf{m}_{y,T}\left(  \delta\right)
:y\in\mathcal{M}\right\}  \smallskip
\]

\begin{proposition}
\label{p-uq-gsp}Fix $T>0$. Let the assumption $(H_{GSP})$ be satisfied and%
\[
int\left(  Dom\left(  A\right)  \right)  \neq\emptyset.
\]
Then, there exists a positive constant $C_{\mathcal{M}}$ such that$\newline%
\left(  a\right)  $\quad If $m\in\mathcal{M}$ and $\left(  x,k\right)
=\mathcal{GSP}\left(  A;x_{0},m\right)  $ then%
\begin{equation}
\left\Vert x\right\Vert _{T}^{2}+\left\updownarrow k\right\updownarrow
_{T}\leq C_{\mathcal{M}}(1+\left\vert x_{0}\right\vert ^{2}). \label{gspa1}%
\end{equation}
$\left(  b\right)  $\quad If $m,\hat{m}\in\mathcal{M}$, $\left(  x,k\right)
=\mathcal{GSP}\left(  A;x_{0},m\right)  $ and $(\hat{x},\hat{k})=\mathcal{GSP}%
(A;\hat{x}_{0},\hat{m})$ then%
\begin{equation}
\left\Vert x-\hat{x}\right\Vert _{T}\leq C_{\mathcal{M}}\left(  1+\left\vert
x_{0}\right\vert +\left\vert \hat{x}_{0}\right\vert \right)  (\left\vert
x_{0}-\hat{x}_{0}\right\vert +\Vert m-\hat{m}\Vert_{T}^{1/2}). \label{gsp2b}%
\end{equation}
In particular, the uniqueness follows, that is, if $x_{0}=\hat{x}_{0}$ and
$m=\hat{m}$ then $\left(  x,k\right)  =(\hat{x},\hat{k}).$
\end{proposition}

\begin{proof}
$\left(  a\right)  \;$In the sequel we fix arbitrary $u_{0}\in\mathbb{H}$ and
$0<r_{0}\leq1$ such that%
\[
\bar{B}\left(  u_{0},r_{0}\right)  \subset Dom\left(  A\right)
\]
and%
\[
A_{u_{0},r_{0}}^{\#}:=\sup\left\{  \left\vert \hat{u}\right\vert :\hat{u}\in
A\left(  u_{0}+r_{0}v\right)  ,\;\left\vert v\right\vert \leq1\right\}
<\infty.
\]
If in (\ref{integr def sol GSP}) we consider $z=u_{0}+r_{0}v$, $\left\vert
v\right\vert \leq1$ and $z^{\ast}\in A\left(  z\right)  $, then $\left\vert
z^{\ast}\right\vert \leq A_{u_{0},r_{0}}^{\#}$ and we infer%
\begin{equation}
r_{0}d\left\updownarrow k\right\updownarrow _{t}\leq\left\langle x\left(
t\right)  -u_{0},dk\left(  t\right)  \right\rangle +A_{u_{0},r_{0}}%
^{\#}\left[  r_{0}+\left\vert x\left(  t\right)  -u_{0}\right\vert \right]
dt.\label{adiez}%
\end{equation}
Let $\delta_{0}=\delta_{0,\mathcal{M}}>0$ be defined by%
\[
\delta_{0}+\mathbf{m}_{\mathcal{M},T}\left(  \delta_{0}\right)  =\frac{r_{0}%
}{4}.
\]

By Energy Equality%
\[
\left\vert x\left(  t\right)  -m\left(  t\right)  -u_{0}\right\vert ^{2}+2%
{\displaystyle\int_{0}^{t}}
\left\langle x\left(  r\right)  -u_{0},dk\left(  r\right)  \right\rangle
=\left\vert x_{0}-u_{0}\right\vert ^{2}+2%
{\displaystyle\int_{0}^{t}}
\left\langle m\left(  r\right)  ,dk\left(  r\right)  \right\rangle
\]
and, using (\ref{adiez}), we obtain%
\[
\left\vert x\left(  t\right)  -m\left(  t\right)  -u_{0}\right\vert
^{2}+2r_{0}\left\updownarrow k\right\updownarrow _{t}\leq\left\vert
x_{0}-u_{0}\right\vert ^{2}+2%
{\displaystyle\int_{0}^{t}}
\left\langle m\left(  r\right)  ,dk\left(  r\right)  \right\rangle
+2A_{u_{0},r_{0}}^{\#}%
{\displaystyle\int_{0}^{t}}
\left[  r_{0}+\left\vert x\left(  r\right)  -u_{0}\right\vert \right]  dr.
\]

Let $n_{0}=\left\lceil \frac{T}{\delta_{0}}\right\rceil $ and consider the
partition $0=t_{0}<t_{1}<...<t_{n_{0}}=t,\;t_{i+1}-t_{i}=\dfrac{t}{n_{0}}%
\leq\delta_{0},\;i=\overline{0,n_{0}-1}$ ($\left\lceil a\right\rceil $ is the
smallest integer greater or equal to $a\in\mathbb{R}$). Then%
\[%
\begin{array}
[c]{l}%
{\displaystyle\int_{0}^{t}}
\left\langle m\left(  r\right)  ,dk\left(  r\right)  \right\rangle =%
{\displaystyle\sum\limits_{i=0}^{n_{0}-1}}
{\displaystyle\int\nolimits_{t_{i}}^{t_{i+1}}}
\left\langle m\left(  r\right)  -m\left(  t_{i}\right)  ,dk\left(  r\right)
\right\rangle +%
{\displaystyle\sum\limits_{i=0}^{n_{0}-1}}
\langle m\left(  t_{i}\right)  ,k\left(  t_{i+1}\right)  -k\left(
t_{i}\right)  \rangle\\
\quad\quad\quad\leq\mathbf{m}_{\mathcal{M},T}(\delta_{0})\left\updownarrow
k\right\updownarrow _{t}+%
{\displaystyle\sum\limits_{i=0}^{n_{0}-1}}
\langle m\left(  t_{i}\right)  ,m\left(  t_{i+1}\right)  -x\left(
t_{i+1}\right)  +u_{0}-m\left(  t_{i}\right)  +x\left(  t_{i}\right)
-u_{0}\rangle\medskip\\
\quad\quad\quad\leq\dfrac{r_{0}}{4}\left\updownarrow k\right\updownarrow
_{t}+2(n_{0}+1)\left\Vert m\right\Vert _{t}\left\Vert x-u_{0}-m\right\Vert
_{t}~.
\end{array}
\]
Hence%
\begin{align*}
\left\vert x\left(  t\right)  -m\left(  t\right)  -u_{0}\right\vert
^{2}+\dfrac{3r_{0}}{2}\left\updownarrow k\right\updownarrow _{t}  &
\leq\left\vert x_{0}-u_{0}\right\vert ^{2}+\left[  4\left(  n_{0}+1\right)
\left\Vert m\right\Vert _{t}+2tA_{u_{0},r_{0}}^{\#}\right]  \left\Vert
x-u_{0}-m\right\Vert _{t}\\
&  +2\left(  t+t\left\Vert m\right\Vert _{t}\right)  A_{u_{0},r_{0}}^{\#},
\end{align*}
which implies (\ref{gspa1}), where $C_{\mathcal{M}}=C(T,u_{0},r_{0}%
,A_{u_{0},r_{0}}^{\#},\delta_{0},\left\Vert \mathcal{M}\right\Vert _{T})$.

\noindent$\left(  b\right)  \;$By ordinary differential calculus and
(\ref{gspa1}) we infer%
\[%
\begin{array}
[c]{l}%
|x\left(  t\right)  -m\left(  t\right)  -\hat{x}\left(  t\right)  +\hat
{m}\left(  t\right)  |^{2}+2%
{\displaystyle\int_{0}^{t}}
\langle x\left(  r\right)  -\hat{x}\left(  r\right)  ,dk\left(  r\right)
-d\hat{k}\left(  r\right)  \rangle\\
\quad\quad\quad\quad\quad\quad\quad\quad\quad\quad\quad=\left\vert x_{0}%
-\hat{x}_{0}\right\vert ^{2}+2%
{\displaystyle\int_{0}^{t}}
\langle m\left(  r\right)  -\hat{m}\left(  r\right)  ,dk\left(  r\right)
-d\hat{k}\left(  r\right)  \rangle\smallskip\\
\quad\quad\quad\quad\quad\quad\quad\quad\quad\quad\quad\leq\left\vert
x_{0}-\hat{x}_{0}\right\vert ^{2}+2\Vert m-\hat{m}\Vert_{T}[\updownarrow
\!k\!\updownarrow_{T}+\updownarrow\!\hat{k}\!\updownarrow_{T}]\medskip\\
\quad\quad\quad\quad\quad\quad\quad\quad\quad\quad\quad\leq\left\vert
x_{0}-\hat{x}_{0}\right\vert ^{2}+4C_{\mathcal{M}}\Vert m-\hat{m}\Vert
_{T}(1+\left\vert x_{0}\right\vert ^{2}+\left\vert \hat{x}_{0}\right\vert
^{2}).
\end{array}
\]
On the other hand,%
\begin{align*}
|x\left(  t\right)  -m\left(  t\right)  -\hat{x}\left(  t\right)  +\hat
{m}\left(  t\right)  |^{2}  &  \geq\frac{1}{2}\left\vert x\left(  t\right)
-\hat{x}\left(  t\right)  \right\vert ^{2}-\Vert m-\hat{m}\Vert_{T}^{2}\\
&  \geq\frac{1}{2}\left\vert x\left(  t\right)  -\hat{x}\left(  t\right)
\right\vert ^{2}-2\left\Vert \mathcal{M}\right\Vert _{T}\Vert m-\hat{m}%
\Vert_{T}%
\end{align*}
Combining these last two inequalities with (\ref{pozitivitate max mon}), we
deduce%
\[
|x\left(  t\right)  -\hat{x}\left(  t\right)  |^{2}\leq2\left\vert x_{0}%
-\hat{x}_{0}\right\vert ^{2}+4\left\Vert \mathcal{M}\right\Vert _{T}\Vert
m-\hat{m}\Vert_{T}+8C_{\mathcal{M}}\Vert m-\hat{m}\Vert_{T}(1+\left\vert
x_{0}\right\vert ^{2}+\left\vert \hat{x}_{0}\right\vert ^{2})
\]
and (\ref{gsp2b}) easily follows, with a constant $\hat{C}_{\mathcal{M}}$; the
two relations (\ref{gspa1}) and (\ref{gsp2b}) can be written with a common
constant $C_{\mathcal{M}}:=\max\{C_{\mathcal{M}},\hat{C}_{\mathcal{M}}%
\}$.\hfill
\end{proof}

\begin{theorem}
\label{Existence for GSP}Under the assumptions $(H_{GSP})$, if we have also
$int\left(  Dom\left(  A\right)  \right)  \neq\emptyset$, then the generalized
convex Skorohod problem (\ref{3}) has a unique solution $\left(  x,k\right)  $
and estimates (\ref{gspa1}) and (\ref{gsp2b}) hold.
\end{theorem}

\begin{proof}
The uniqueness and estimates (\ref{gspa1}) and (\ref{gsp2b}) have been
obtained in the above result. It suffices to prove the existence on an
arbitrary fixed interval $\left[  0,T\right]  $.

Let $x_{0,n}\in Dom(A)$ and $m_{n}\in C^{\infty}\left(  \left[  0,T\right]
;\mathbb{H}\right)  $ be such that%
\[
x_{0,n}\rightarrow x_{0}\quad\text{in \ }\mathbb{H}\quad\quad\text{and}%
\quad\text{\quad}m_{n}\rightarrow m\quad\text{in\ \ }C\left(  \left[
0,T\right]  ;\mathbb{H}\right)  .
\]
Notice that $\mathcal{M}=\left\{  m,m_{1},m_{2},\ldots\right\}  $ is a bounded
equicontinuous subset of $C\left(  \left[  0,T\right]  ;\mathbb{H}\right)  $.
We set $\alpha(\varepsilon)=\mathbf{m}_{\mathcal{M},T}\left(  \varepsilon
\right)  $ and let $\hat{J}$ (resp. $\hat{J}_{n}$)$:\mathbb{H}\times
\mathbb{X}\times\mathbb{Y}_{R}\times\mathbb{X}\rightarrow]-\infty,+\infty]$ be
the functions defined by (\ref{fcJ}) associated to $\left(  x_{0},m,A\right)
$ (and resp. $\left(  x_{0,n},m_{n},A\right)  $). Then%
\begin{align*}
\hat{J}\left(  a,x,k,\mu\right)   &  =\hat{J}_{n}\left(  a,x,k,\mu\right)
-\left\vert a-x_{0,n}\right\vert ^{2}-2R\left\Vert \mu-m_{n}\right\Vert
_{T}+\left\vert a-x_{0}\right\vert ^{2}\smallskip\\
&  +\sup_{\nu\in C_{\alpha}}\left\{  \left\langle \!\left\langle \mu
-\nu,k\right\rangle \!\right\rangle -R\left\Vert \nu-m\right\Vert
_{T}\right\}  -\sup_{\nu\in C_{\alpha}}\left\{  \left\langle \!\left\langle
\mu-\nu,k\right\rangle \!\right\rangle -R\left\Vert \nu-m_{n}\right\Vert
_{T}\right\}  \smallskip\\
&  \leq\hat{J}_{n}\left(  a,x,k,\mu\right)  -\left\vert a-x_{0,n}\right\vert
^{2}-2R\left\Vert \mu-m_{n}\right\Vert _{T}+\left\vert a-x_{0}\right\vert
^{2}\smallskip\\
&  +R\sup_{\nu\in C_{\alpha}}\left\{  \left\Vert \nu-m_{n}\right\Vert
_{T}-\left\Vert \nu-m\right\Vert _{T}\right\}  \smallskip\\
&  \leq\hat{J}_{n}\left(  a,x,k,\mu\right)  -\left\vert a-x_{0,n}\right\vert
^{2}-2R\left\Vert \mu-m_{n}\right\Vert _{T}+\left\vert a-x_{0}\right\vert
^{2}+R\left\Vert m-m_{n}\right\Vert _{T}~.
\end{align*}
In particular,%
\begin{equation}
\hat{J}\left(  x_{0,n},x,k,m_{n}\right)  \leq\hat{J}_{n}\left(  x_{0,n}%
,x,k,m_{n}\right)  +\left\vert x_{0,n}-x_{0}\right\vert ^{2}+R\left\Vert
m-m_{n}\right\Vert _{T}~. \label{relation J Jc}%
\end{equation}
By a classical result (see Barbu \cite{B75}, Theorem 2.2) there exist
$x_{n}\in C\left(  \left[  0,T\right]  ;\mathbb{H}\right)  $ and\newline%
$h_{n}\in L^{1}\left(  0,T;\mathbb{H}\right)  $, $h_{n}\left(  t\right)  \in
Ax_{n}(t),$ a.e. $t\in\left[  0,T\right]  $, such that%
\begin{equation}
x_{n}\left(  t\right)  +%
{\displaystyle\int_{0}^{t}}
h_{n}\left(  s\right)  ds=x_{0,n}+m_{n}\left(  t\right)  .
\label{ecuatie aproximanta}%
\end{equation}
If we denote $k_{n}\left(  t\right)  =%
{\displaystyle\int_{0}^{t}}
h_{n}\left(  s\right)  ds$, then $\left(  x_{n},k_{n}\right)  \in\mathcal{A}$
and therefore, by Fitzpatrick's Theorem, $\mathcal{H}\left(  x_{n}%
,k_{n}\right)  =\left\langle \!\left\langle x_{n},k_{n}\right\rangle
\!\right\rangle .$

\noindent Then, using Proposition \ref{p-uq-gsp}, there exists a positive
constant $\mathcal{C}$, not depending on $n$, such that, for all
$n,j\in\mathbb{N}^{\ast}$,%
\begin{align*}
\left\Vert x_{n}\right\Vert _{T}^{2}+\left\updownarrow k_{n}\right\updownarrow
_{T}  &  \leq\mathcal{C}\text{ and}\\
\left\Vert x_{n}-x_{j}\right\Vert _{T}  &  \leq\mathcal{C}(\left\vert
x_{0,n}-x_{0,j}\right\vert +\Vert m_{n}-m_{j}\Vert_{T}^{1/2}).
\end{align*}
Hence, there exists $x\in C\left(  \left[  0,T\right]  ;\mathbb{H}\right)  $
such that, as $n\rightarrow\infty$,%
\[
x_{n}\rightarrow x\quad\text{in\ \ }C(\left[  0,T\right]  ;\overline
{Dom(A)}).
\]
Let%
\[
k\left(  t\right)  =x_{0}+m\left(  t\right)  -x\left(  t\right)  .
\]
We deduce that%
\[
k_{n}=x_{0,n}+m_{n}-x_{n}\;\longrightarrow\;k\quad\text{in }C\left(  \left[
0,T\right]  ;\mathbb{H}\right)
\]
and clearly follows%
\[
k\in BV\left(  \left[  0,T\right]  ;\mathbb{H}\right)  ,\quad\left\updownarrow
k\right\updownarrow _{T}\leq\mathcal{C}.
\]
Setting $R=\mathcal{C}$, the quantities $\hat{J}\left(  x_{0,n},x_{n}%
,k_{n},m_{n}\right)  $ and $\hat{J}_{n}\left(  x_{0,n},x_{n},k_{n}%
,m_{n}\right)  $ are well defined. Moreover, by Proposition \ref{p1}, $\hat
{J}_{n}\left(  x_{0,n},x_{n},k_{n},m_{n}\right)  =0$. Passing to
$\liminf\limits_{n\rightarrow+\infty}$ in (\ref{relation J Jc}), the
lower-semicontinuity of $\hat{J}$ implies%
\[
0\leq\hat{J}\left(  x_{0},x,k,m\right)  \leq\liminf\limits_{n\rightarrow
+\infty}\hat{J}\left(  x_{0,n},x_{n},k_{n},m_{n}\right)  =0,
\]
that is, there exists a minimum point for which $\hat{J}$ is zero. By
Proposition \ref{p1} (-$\left(  c\right)  $) we infer that the generalized
convex Skorohod problem (\ref{3}) has a solution.\hfill
\end{proof}

\begin{remark}
We highlight that the existence problem is reduced to the minimization of a
specific l.s.c. convex function on a bounded closed convex subset of
$\mathbb{H}\times\mathbb{X}\times BV\left(  \left[  0,T\right]  ;\mathbb{H}%
\right)  \times\mathbb{X}$. Indeed, via Proposition \ref{p1} (-$\left(
c\right)  $), the minimization of $\hat{J}$ is on the set $\mathbb{H}%
_{\rho_{0}}\times\mathbb{X}_{R}\times\mathbb{Y}_{R}\times C_{\alpha}$, where%
\[
\mathbb{H}_{\rho_{0}}=\left\{  h\in\mathbb{H}:\left\vert h\right\vert \leq
\rho_{0}:=\sup\{|x_{0}|,|x_{0,n}|:n\in\mathbb{N}^{\ast}\}\right\}  ,
\]
$\mathbb{X}_{R}=\{x\in\mathbb{X}:\left\Vert x\right\Vert _{T}\leq R\}$ and
$R=\mathcal{C}$. Classical results (see Zeidler \cite{Z}, Theorem 38.A)
establish sufficient conditions for a functional defined on a subset of a
reflexive Banach space to attain its minimum.
\end{remark}

We note that, in the framework of Hilbert spaces, the assumption
$int(Dom(A))\neq\emptyset$ from the above results is fairly restrictive. One
can renounce at this condition, but we have to consider a stronger assumption
on $m$ and, moreover, to weaken the notion of solution for the generalized
Skorohod problem (\ref{3}). Therefore, along $\mathbb{H}$, we consider
$\left(  \mathbb{V},\left\Vert \cdot\right\Vert _{\mathbb{V}}\right)  $ a real
separable Banach space with separable dual $\left(  \mathbb{V}^{\ast
},\left\Vert \cdot\right\Vert _{\mathbb{V}^{\ast}}\right)  $ such that%
\[
\mathbb{V}\subset\mathbb{H}\cong\mathbb{H}^{\ast}\subset\mathbb{V}^{\ast},
\]
where the embeddings are continuous, with dense range (the duality paring
$\left(  \mathbb{V}^{\ast},\mathbb{V}\right)  $ is denoted also by
$\left\langle \cdot,\cdot\right\rangle $, and, for $k:[0,\infty
)\longrightarrow\mathbb{V}^{\ast}$, $k\left(  0\right)  =0$, we use the
adequate notation ${\large \updownarrow}\!\!{\large \updownarrow
}k{\large \updownarrow}\!\!{\large \updownarrow}_{\ast T}=\left\Vert
k\right\Vert _{BV\left(  \left[  0,T\right]  ;\mathbb{V}^{\ast}\right)  }$).

Reconsider the multivalued monotone differential equation (\ref{3}) under the
assumptions%
\[
\bar{H}_{GSP}:\quad\left\{
\begin{array}
[c]{l}%
H_{GSP}:\left(  i\right)  \quad\text{and}\quad\left(  ii\right)  ,\medskip\\
\left(  iii^{\prime}\right)  \quad m:[0,\infty)\longrightarrow\mathbb{V}\text{
is continuous and }m\left(  0\right)  =0\text{.}%
\end{array}
\right.
\]

\begin{definition}
\label{extended definition}A continuous function $x:[0,\infty)\rightarrow
\mathbb{H}$ is a solution of Eq.(\ref{3}) if

\begin{itemize}
\item[$\left(  i\right)  $] there exist the sequences $\{x_{0,n}\}\subset
Dom(A)$ and $m_{n}:[0,\infty)\longrightarrow\mathbb{V},$ $m_{n}\left(
0\right)  =0$ of $C^{1}-$continuous functions satisfying, for all $T>0$,%
\[
\left\vert x_{0,n}-x_{0}\right\vert +\left\Vert m_{n}-m\right\Vert _{C\left(
\left[  0,T\right]  ;\mathbb{V}\right)  }\rightarrow0,\ \text{as }%
n\rightarrow\infty\text{,}%
\]

\item[$\left(  ii\right)  $] there exist $x_{n}\in C([0,\infty);\overline
{Dom(A)}\mathbb{)},$ $k_{n}\in C([0,\infty);\mathbb{H)\cap}BV_{0,loc}\left(
\mathbb{R}_{+};\mathbb{V}^{\ast}\right)  ,$ $k_{n}\left(  0\right)  =0,$ and a
function $k$ such that
\[
x_{n}\left(  t\right)  +k_{n}\left(  t\right)  =x_{0,n}+m_{n}\left(  t\right)
,\ \forall t\geq0
\]
and, for all $T>0$,%
\[%
\begin{array}
[c]{ll}%
\left(  a\right)  \quad & \left\Vert x_{n}-x\right\Vert _{T}+\left\Vert
k_{n}-k\right\Vert _{T}\rightarrow0,\ \text{as }n\rightarrow\infty
,\smallskip\\
\left(  b\right)  \quad & \sup\limits_{n\in\mathbb{N}^{\ast}}%
{\large \updownarrow}\!\!{\large \updownarrow}k_{n}{\large \updownarrow
}\!\!{\large \updownarrow}_{\ast T}<\infty,\smallskip\\
\left(  c\right)  \quad &
{\displaystyle\int_{s}^{t}}
\left\langle x_{n}\left(  r\right)  -z,dk_{n}\left(  r\right)  -z^{\ast
}dr\right\rangle \geq0,\ \forall~\left(  z,z^{\ast}\right)  \in A,\ \forall
0\leq s\leq t\leq T\text{.}%
\end{array}
\]

\end{itemize}

\noindent(Without confusion, the uniqueness of $k$ will permit us to call the
pair $\left(  x,k\right)  $ solution of the generalized Skorohod problem
(\ref{3}) and we write $\left(  x,k\right)  =\mathcal{GSP}\left(
A;x_{0},m\right)  $.)
\end{definition}

\begin{remark}
If $\left(  x,k\right)  =\mathcal{GSP}\left(  A;x_{0},m\right)  $ then we
clearly have

\begin{itemize}
\item[$\left(  iii\right)  $] $x\left(  t\right)  \in\overline{Dom(A)}$, for
all $t\geq0,$

\item[$\left(  iv\right)  $] $k\in C([0,\infty);\mathbb{H)\cap}BV_{0,loc}%
\left(  \mathbb{R}_{+};\mathbb{V}^{\ast}\right)  $, $k\left(  0\right)  =0$ and

\item[$\left(  v\right)  $] $x\left(  t\right)  +k\left(  t\right)
=x_{0}+m\left(  t\right)  ,\ \forall t\geq0$.
\end{itemize}
\end{remark}

Replacing now the condition $int(Dom(A))\neq\emptyset$ we obtain (see, for
example, R\u{a}\c{s}canu \cite{R}, Theorem 2.3) the following result of
existence and uniqueness of a solution for the generalized Skorohod problem
(\ref{3}).

\begin{theorem}
\label{Existence in VHV*}Under the hypothesis $\left(  \bar{H}_{GSP}\right)
$, if there exist $h_{0}\in\mathbb{H}$ and $r_{0},a_{1},a_{2}>0$ such that%
\begin{equation}
r_{0}\left\Vert z^{\ast}\right\Vert _{\mathbb{V}^{\ast}}\leq\left\langle
z^{\ast},z-h_{0}\right\rangle +a_{1}\left\vert z\right\vert ^{2}+a_{2}%
,\quad\forall\left(  z,z^{\ast}\right)  \in A \label{condition for A}%
\end{equation}
then the differential equation (\ref{3}) has a unique solution $\left(
x,k\right)  $ in the sense of Definition \ref{extended definition}. Moreover,
for all $T>0$,

\begin{itemize}
\item[$\left(  a\right)  $] if $\left(  x,k\right)  =\mathcal{GSP}\left(
A;x_{0},m\right)  $ and $(\hat{x},\hat{k})=\mathcal{GSP}\left(  A;\hat{x}%
_{0},\hat{m}\right)  $, then there exists a positive constant $C$ such that%
\[
\left\Vert x-\hat{x}\right\Vert _{T}^{2}\leq C\left[  \left\vert x_{0}-\hat
{x}_{0}\right\vert ^{2}+\left\Vert m-\hat{m}\right\Vert _{T}^{2}+\left\Vert
m-\hat{m}\right\Vert _{C\left(  \left[  0,T\right]  ;\mathbb{V}\right)
}{\large \updownarrow}\!\!{\large \updownarrow k-\hat{k}\updownarrow
}\!\!{\large \updownarrow}_{\ast T}\right]  \text{ and}%
\]

\item[$\left(  b\right)  $] for every equiuniform continuous subset
$\mathcal{M}\subset C\left(  \left[  0,T\right]  ;\mathbb{V}\right)  ,$
$m\in\mathcal{M},$ there exists $C_{0}=C_{0}\left(  r_{0},h_{0},a_{1}%
,a_{2},T,\mathcal{N}_{\mathcal{M}}\right)  >0$ for which%
\[
\left\Vert x\right\Vert _{T}^{2}+{\large \updownarrow}\!\!{\large \updownarrow
k\updownarrow}\!\!{\large \updownarrow}_{\ast T}\leq C_{0}\left[  1+\left\vert
x_{0}\right\vert ^{2}+\left\Vert m\right\Vert _{T}^{2}\right]  .
\]
(Here $\mathcal{N}_{\mathcal{M}}$ is the constant of equiuniform continuity
given by $\sup\{\left\Vert f\left(  t\right)  -f\left(  s\right)  \right\Vert
_{\mathbb{V}}:\left\vert t-s\right\vert \leq T/\mathcal{N}_{\mathcal{M}}\}\leq
r_{0}/4,\ \forall f\in\mathcal{M}.$)
\end{itemize}
\end{theorem}

From R\u{a}\c{s}canu \cite{R} we mention three situations when the relation
(\ref{condition for A}) is satisfied:

\begin{itemize}
\item[$\left(  a\right)  $] $A=A_{0}+\partial\varphi,\mathit{\ }$where\textit{
}$A_{0}:\mathbb{H}\rightarrow\mathbb{H}$ is a continuous monotone operator on%
\
$\mathbb{H}$ and $\varphi:\mathbb{H}\rightarrow]-\infty,+\infty]$ is a proper
convex l.s.c. function for which there exist $h_{0}\in\mathbb{H}%
,\;R_{0}>0,\;a_{0}>0$ such that%
\[
\varphi\left(  h_{0}+x\right)  \leq a_{0},\;\forall x\in\mathbb{V}%
,\;\left\Vert x\right\Vert _{\mathbb{V}}\leq R_{0}.
\]

\item[$\left(  b\right)  $] 

\begin{itemize}
\item[$\circ$] There exists$\,$a separable Banach space $\mathbb{U}$ such that%
\
$\mathbb{U}\subset\mathbb{H}\subset\mathbb{U}^{\ast}$ densely and continuously
and%
\
$\mathbb{U}\cap\mathbb{V}$%
\
is dense in%
\
$\mathbb{V}$,

\item[$\circ$] $A:\mathbb{H}\rightrightarrows\mathbb{H}$ is a maximal monotone
operator with%
\
$Dom(A)\subset\mathbb{U}$,

\item[$\circ$] $\exists a,\lambda\in\mathbb{R},$\textit{\ }$a>0$, such that
for all $\left(  x_{1},y_{1}\right)  ,\ \left(  x_{2},y_{2}\right)  \in A$%
\[
(y_{1}-y_{2},x_{1}-x_{2})+\lambda\left\vert x_{1}-x_{2}\right\vert ^{2}\geq
a\left\Vert x_{1}-x_{2}\right\Vert _{\mathbb{V}}^{2},
\]

\item[$\circ$] $\exists h_{0}\in\mathbb{U},\;\exists r_{0},a_{0}>0$%
\
such that%
\[
h_{0}+r_{0}e\in Dom(A)\quad\text{and}\quad\left\Vert A^{0}\left(  h_{0}%
+r_{0}e\right)  \right\Vert _{\mathbb{U}^{\ast}}\leq r_{0},
\]%
\
for%
\textit{ }%
all%
\
$e\in\mathbb{U}\cap\mathbb{V},$ $\left\Vert e\right\Vert _{\mathbb{V}}=1$,
where $A^{0}x:=\Pr_{Ax}0$.
\end{itemize}

\item[$\left(  c\right)  $] $A$%
\
is a maximal monotone with%
\
$int(Dom(A))\neq\emptyset$ and $\mathbb{V=H}$.
\end{itemize}

\subsection{Maximal monotone SDE with additive noise}

Consider now the following stochastic differential equation (for short SDE),
where by $B$ we denote the $\mathbb{H}_{0}$-Wiener process defined in Section
\ref{sfram},%
\begin{equation}
\left\{
\begin{array}
[c]{l}%
dX_{t}+AX_{t}(dt)\ni G_{t}dB_{t}~,\medskip\\
X_{0}=\xi,\quad t\in\left[  0,T\right]  ,
\end{array}
\right.  \label{SDE-AN}%
\end{equation}
where%
\[
\left(  H_{MSDE}\right)  :\quad\left\{
\begin{array}
[c]{rl}%
\left(  i\right)  & A:\mathbb{H}\rightrightarrows\mathbb{H}\text{ is a maximal
monotone operator,}\medskip\\
\left(  ii\right)  & \xi\in L^{0}(\Omega,\mathcal{F}_{0},\mathbb{P}%
;\overline{Dom(A)}),\medskip\\
\left(  iii\right)  & G\in\Lambda_{\mathbb{H}\times\mathbb{H}_{0}}^{2}~.
\end{array}
\right.
\]
Setting $\mathbb{X}=L^{2}\left(  \Omega;C\left(  \left[  0,T\right]
;\mathbb{H}\right)  \right)  $, the space $L^{2}\left(  \Omega;BV_{0}\left(
\left[  0,T\right]  ;\mathbb{H}\right)  \right)  $ is a linear subspace of the
dual of $\mathbb{X}$ and, the natural duality%
\[
\left(  X,K\right)  \mapsto\mathbb{E}%
{\displaystyle\int_{0}^{T}}
\left\langle X_{t},dK_{t}\right\rangle
\]
between these two suggests to use the notation $\mathbb{X}^{\ast}$ for
$L^{2}\left(  \Omega;BV_{0}\left(  \left[  0,T\right]  ;\mathbb{H}\right)
\right)  $, even it is not the entire dual space. On $\mathbb{X}$ we shall
consider the strong topology and on $\mathbb{X}^{\ast}$ the $w^{\ast}%
$-topology. Let $\mathcal{A}$ the realization of $A$ on $\mathbb{X}%
\times\mathbb{X}^{\ast}.$

\begin{definition}
\label{definition with H}By a solution of Eq.(\ref{SDE-AN}) we understand a
pair of stochastic processes%
\[
\left(  X,K\right)  \in L^{0}\left(  \Omega;C\left(  \left[  0,T\right]
;\mathbb{H}\right)  \right)  \times\left[  L^{0}\left(  \Omega;C\left(
\left[  0,T\right]  ;\mathbb{H}\right)  \right)  \cap L^{0}\left(
\Omega;BV_{0}\left(  \left[  0,T\right]  ;\mathbb{H}\right)  \right)  \right]
,
\]
satisfying, $\mathbb{P}$-a.s. $\omega\in\Omega$, for all $0\leq s\leq t\leq
T$,%
\[%
\begin{array}
[c]{ll}%
\left(  c_{1}\right)  \quad & X_{t}\in\overline{Dom(A)},\\
\left(  c_{2}\right)  \quad & X_{t}+K_{t}=\xi+%
{\displaystyle\int_{0}^{t}}
G_{s}dB_{s}\text{ and}\\
\left(  c_{3}\right)  \quad &
{\displaystyle\int_{s}^{t}}
\left\langle X_{r}-u,dK_{r}-vdr\right\rangle \geq0,\text{ }\forall(u,v)\in A.
\end{array}
\]

\end{definition}

\noindent Clearly,%
\[
\left(  X(\omega,\cdot),K(\omega,\cdot)\right)  =\mathcal{GSP}\left(
A;\xi(\omega),M(\omega,\cdot)\right)  ,\text{ }\mathbb{P}\text{-}%
a.s.\ \omega\in\Omega,
\]
where $M_{t}=%
{\displaystyle\int_{0}^{t}}
G_{s}dB_{s}\in\mathcal{M}^{2}(0,T;\mathbb{H})$. Consequently, under the
hypothesis $\left(  H_{MSDE}\right)  $, if $int(Dom(A))\neq\emptyset$ then by
Theorem \ref{Existence for GSP} there exists a unique solution $\left(
X,K\right)  $ (in the sense of Definition \ref{definition with H}) for
Eq.(\ref{SDE-AN}). Moreover, if%

\[
\mathbb{E}\left\vert \xi\right\vert ^{4}+\mathbb{E}\left(
{\displaystyle\int_{0}^{T}}
\left\Vert G_{t}\right\Vert _{HS}^{2}dt\right)  ^{2}<+\infty
\]
then $X\in L^{4}\left(  \Omega;C\left(  \left[  0,T\right]  ;\mathbb{H}%
\right)  \right)  \subset\mathbb{X}$ and $K\in\mathbb{X\cap X}^{\ast}$ (see
for example Pardoux \& R\u{a}\c{s}canu \cite{PRb}, Proposition 4.22).

In the sequel we define a convex functional whose minimum point coincide with
the solution of Eq.(\ref{SDE-AN}).

Let%
\[
\mathbb{S}=L^{2}\left(  \Omega,\mathcal{F}_{0},\mathbb{P};\mathbb{H}\right)
\times\mathbb{X}\times\mathbb{X}^{\ast}\times\Lambda_{\mathbb{H}%
\times\mathbb{H}_{0}}^{2}~.
\]

Define, for each $\left(  U,U^{\ast}\right)  \in\mathcal{A}$,%
\[
J_{\left(  U,U^{\ast}\right)  }:\mathbb{S}\rightarrow\mathbb{R}%
\]
by%
\[%
\begin{array}
[c]{ll}%
J_{\left(  U,U^{\ast}\right)  }\left(  \eta,X,K,g\right)  & =\dfrac{1}%
{2}\mathbb{E}\left\vert \eta-\xi\right\vert ^{2}+\dfrac{1}{2}\mathbb{E}%
{\displaystyle\int_{0}^{T}}
\left\Vert g_{t}-G_{t}\right\Vert _{HS}^{2}dt\smallskip\\
& +\mathbb{E}%
{\displaystyle\int_{0}^{T}}
\left[  \left\langle U_{t},dK_{t}\right\rangle +\left\langle X_{t}%
,dU_{t}^{\ast}\right\rangle -\left\langle U_{t},dU_{t}^{\ast}\right\rangle
-\left\langle X_{t},dK_{t}\right\rangle \right]
\end{array}
\]
and $\hat{J}:\mathbb{S}\rightarrow]-\infty,+\infty]$%
\begin{align*}
\hat{J}\left(  \eta,X,K,g\right)   &  =\sup_{\left(  U,U^{\ast}\right)
\in\mathcal{A}}J_{\left(  U,U^{\ast}\right)  }\left(  \eta,X,K,g\right) \\
&  =\dfrac{1}{2}\mathbb{E}\left\vert \eta-\xi\right\vert ^{2}+\mathcal{H}%
\left(  X,K\right)  -\ll X,K\gg+\dfrac{1}{2}\mathbb{E}%
{\displaystyle\int_{0}^{T}}
\left\Vert g_{t}-G_{t}\right\Vert _{HS}^{2}dt,
\end{align*}
where $\mathcal{H}:\mathbb{X}\times\mathbb{X}^{\ast}\rightarrow]-\infty
,+\infty]$ is the Fitzpatrick function associated to the maximal monotone
operator $\mathcal{A}$. It is clear that

\begin{remark}
$\hat{J}:\mathbb{S}\rightarrow]-\infty,+\infty]$ is a lower semicontinuous
function as supremum of continuous functions.
\end{remark}

Since $\mathcal{H}\left(  X,K\right)  \geq\left\langle \!\left\langle
X,K\right\rangle \!\right\rangle $, then we easily deduce

\begin{proposition}
$\hat{J}$ has the following properties:

\begin{itemize}
\item[$\left(  a\right)  $] $\;\hat{J}\left(  \eta,X,K,g\right)  \geq0$, for
all $\left(  \eta,X,K,g\right)  \in\mathbb{S}$.

\item[$\left(  b\right)  $] $\;\hat{J}\left(  \eta,X,K,g\right)  =0$ iff
$\eta=\xi,$ $g=G$ and $K\in\mathcal{A}\left(  X\right)  $.

\item[$\left(  c\right)  $] $\;$Let $R>0$. The restriction of $\hat{J}$ to the
bounded closed convex set%
\[%
\begin{array}
[c]{ccc}%
\mathcal{L} & = & \left\{  \left(  \eta,X,K,g\right)  \in\mathbb{S}%
:X_{t}+K_{t}=\eta+%
{\displaystyle\int_{0}^{t}}
g_{s}dB_{s},\ \forall t\in\left[  0,T\right]  ,\right.  \medskip\\
&  & \left.  \mathbb{E}\left\vert \eta\right\vert ^{2}+\mathbb{E}\left\Vert
X\right\Vert _{\mathbb{X}}^{2}+\mathbb{E}{\large \updownarrow}%
K{\large \updownarrow}_{\mathbb{X}^{\ast}}+\mathbb{E}%
{\displaystyle\int_{0}^{T}}
\left\Vert g_{s}\right\Vert _{HS}^{2}ds\leq R\right\}
\end{array}
\]
is a convex l.s.c. function and $\hat{J}\left(  \eta,X,K,g\right)  =0$ iff
$\eta=\xi,$ $g=G$ and $\left(  X,K\right)  $ is the solution of the SDE
(\ref{SDE-AN}).
\end{itemize}
\end{proposition}

\begin{proof}
The points $\left(  a\right)  $ and $\left(  b\right)  $ clearly are
consequences of the properties of the Fitzpatrick function $\mathcal{H}.$ Let
us prove $\left(  c\right)  .$ Since, by Energy Equality%
\[
\frac{1}{2}\mathbb{E}\left\vert X_{T}\right\vert ^{2}+\mathbb{E}%
{\displaystyle\int_{0}^{T}}
\left\langle X_{t},dK_{t}\right\rangle =\frac{1}{2}\mathbb{E}\left\vert
\eta\right\vert ^{2}+\frac{1}{2}\mathbb{E}%
{\displaystyle\int_{0}^{T}}
\left\Vert g_{t}\right\Vert _{HS}^{2}dt
\]
then%
\begin{align*}
\hat{J}\left(  \eta,X,K,g\right)   &  =\dfrac{1}{2}\mathbb{E}\left\vert
\eta-\xi\right\vert ^{2}+\mathcal{H}\left(  X,K\right)  -\left\langle
\!\left\langle X,K\right\rangle \!\right\rangle +\dfrac{1}{2}\mathbb{E}%
{\displaystyle\int_{0}^{T}}
\left\Vert g_{t}-G_{t}\right\Vert _{HS}^{2}dt\\
&  =\dfrac{1}{2}\mathbb{E}\left\vert \xi\right\vert ^{2}-\mathbb{E}%
\left\langle \eta,\xi\right\rangle +\mathcal{H}\left(  X,K\right)  +\dfrac
{1}{2}\mathbb{E}\left\vert X_{T}\right\vert ^{2}\\
&  -\mathbb{E}%
{\displaystyle\int_{0}^{T}}
\left\langle g_{t},G_{t}\right\rangle dt+\dfrac{1}{2}\mathbb{E}%
{\displaystyle\int_{0}^{T}}
\left\Vert G_{t}\right\Vert _{HS}^{2}dt
\end{align*}
and the convexity of $\hat{J}$ on the set $\mathcal{L}$ follows.\hfill\medskip
\end{proof}

To complete this section, we will situate in the extended framework introduced
in the final part of Subsection 2.1. We will consider once again the spaces
$\mathbb{H}$ and $\mathbb{V}$ and we assume that $\mathbb{V}\subset
\mathbb{H}\cong\mathbb{H}^{\ast}\subset\mathbb{V}^{\ast}$, where the
embeddings are continuous with dense range. Concerning the SDE (\ref{SDE-AN}),
the hypothesis $\left(  H_{MSDE}\right)  $ will be replaced by%
\[
\left(  \bar{H}_{MSDE}\right)  :\quad\left\{
\begin{array}
[c]{rl}%
\left(  i\right)  & \left\vert
\begin{array}
[c]{l}%
A:\mathbb{H}\rightrightarrows\mathbb{H}\text{ is a maximal monotone operator
and}\smallskip\\
\text{there exist }h_{0}\in\mathbb{H}\text{ and }r_{0},a_{1},a_{2}>0\text{
such that}\smallskip\\
r_{0}\left\Vert z^{\ast}\right\Vert _{\mathbb{V}^{\ast}}\leq\left\langle
z^{\ast},z-h_{0}\right\rangle +a_{1}\left\vert z\right\vert ^{2}+a_{2},\text{
}\forall\left(  z,z^{\ast}\right)  \in A\medskip
\end{array}
\right. \\
\left(  ii\right)  & \xi\in L^{2}(\Omega,\mathcal{F}_{0},\mathbb{P}%
;\overline{Dom(A)}),\medskip\\
\left(  iii\right)  & G\in\Lambda_{\mathbb{H}\times\mathbb{H}_{0}}^{2}\left(
0,T;\mathcal{L}^{2}\left(  \mathbb{H}_{0},\mathbb{H}\right)  \right)  .
\end{array}
\right.
\]

\begin{definition}
Let $M_{t}:=\int\nolimits_{0}^{t}G_{s}dB_{s}$. A \textit{stochastic process
}$X\in L_{ad}^{0}\left(  \Omega;C\left(  \left[  0,T\right]  ;\mathbb{H}%
\right)  \right)  $ that satisfies, $\mathbb{P}$-a.s., $X_{0}=\xi$ and
$X_{t}\in\overline{Dom(A)},$ $\forall t\in\left[  0,T\right]  $ \textit{is a
(generalized) solution of multivalued SDE }(\ref{SDE-AN}) if there exist%
\[
K\in L_{ad}^{0}\left(  \Omega;C\left(  \left[  0,T\right]  ;\mathbb{H}\right)
\right)  \cap L^{0}\left(  \Omega;BV\left(  0,T;\mathbb{V}^{\ast}\right)
\right)  ,K_{0}=0\mathit{\ }\mathbb{P}\text{-}a.s
\]
and a sequence of stochastic processes $\{M^{n}\}_{n\in\mathbb{N}^{\ast}}$
satisfying%
\begin{equation}
\left\{
\begin{array}
[c]{l}%
M^{n}\in L_{ad}^{2}\left(  \Omega;C\left(  \left[  0,T\right]  ;\mathbb{V}%
\right)  \right)  \cap\mathcal{M}^{2}\left(  0,T;\mathbb{H}\right)
,\smallskip\\
M^{n}\longrightarrow M\text{ in }\mathcal{M}^{2}\left(  0,T;\mathbb{H}\right)
\end{array}
\right.  \label{convergenta}%
\end{equation}
such that, denoting for a.s. $\omega\in\Omega$,%
\[
\left(  X^{n}(\omega,\cdot),K^{n}(\omega,\cdot)\right)  =\mathcal{GSP}\left(
A;\xi\left(  \omega\right)  ,M^{n}\left(  \omega,\cdot\right)  \right)
,\text{ we have}%
\]
$X^{n}\rightarrow X,$ $K^{n}\rightarrow K$ in $L_{ad}^{0}\left(
\Omega,C\left(  \left[  0,T\right]  ;\mathbb{H}\right)  \right)  $ as
$n\rightarrow\infty$ and $\sup\limits_{n}\mathbb{E}{\large \updownarrow
}\!\!{\large \updownarrow}K^{n}{\large \updownarrow}\!\!{\large \updownarrow
}_{\ast T}<+\infty.$

\noindent(Without confusion, the uniqueness of $K$ permits us to call the pair
$(X,K)$ a generalized solution \textit{of the multivalued SDE }(\ref{SDE-AN}).)
\end{definition}

Recall, from R\u{a}\c{s}canu \cite{R}, the following existence result which is
a consequence of the corresponding deterministic case here above.

\begin{theorem}
\textit{Under the assumption }($\bar{H}_{MSDE}$)\textit{\ the problem
}(\ref{SDE-AN}) \textit{has a unique generalized solution }$\left(
X,K\right)  .$ \textit{Moreover the solution satisfies}%
\begin{equation}
\mathbb{E}\sup\limits_{t\in\left[  0,T\right]  }\left\vert X_{t}\right\vert
^{2}+\mathbb{E}\sup\limits_{t\in\left[  0,T\right]  }\left\vert K_{t}%
\right\vert ^{2}+\mathbb{E}\updownarrow\!\left\updownarrow K\right\updownarrow
\!\updownarrow_{\ast T}\leq C_{0}\left[  1+\mathbb{E}\left\vert \xi\right\vert
^{2}+\mathbb{E}%
{\displaystyle\int_{0}^{T}}
\left\Vert G_{t}\right\Vert _{HS}^{2}~dt\right]  \text{,} \label{1 SDE-AN}%
\end{equation}
where $C_{0}=C_{0}\left(  T,r_{0},h_{0},a_{1},a_{2}\right)  >0.$\newline If
$\left(  X,K\right)  $ and $(\tilde{X},\tilde{K})$ are two solutions of
(\ref{SDE-AN}) corresponding to $\left(  \xi,G\right)  $ and, respectively,
$(\tilde{\xi},\tilde{G})$ then%
\begin{equation}
\mathbb{E}\sup\limits_{t\in\left[  0,T\right]  }|X_{t}-\tilde{X}_{t}|^{2}\leq
C\left(  T\right)  \left[  \mathbb{E}|\xi-\tilde{\xi}|^{2}+\mathbb{E}%
{\displaystyle\int_{0}^{T}}
||G_{t}-\tilde{G}_{t}||_{HS}^{2}~dt\right]  . \label{2 SDE-AN}%
\end{equation}

\end{theorem}

\begin{proof}
Since the process $M$ does not have $\mathbb{V}$-valued continuous
trajectories, we use the deterministic result approximating the stochastic
integral by the sequence%
\[
M_{t}^{n}:=\sum_{i=1}^{n}\left\langle M_{t},e_{i}\right\rangle e_{i},
\]
where $\{e_{i};\ i\in\mathbb{N}^{\ast}\}\subset\mathbb{V}$ is an orthonormal
basis in $\mathbb{H}$. By Theorem \ref{Existence in VHV*}, there exists
$(X^{n}\left(  \omega\right)  ,K^{n}\left(  \omega\right)  )=\mathcal{GSP}%
(A;\xi(\omega),M^{n}(\omega)),$ $\mathbb{P}$-$a.s.\;\omega\in\Omega$. It is
not difficult to prove that the following inequalities hold%
\[
\mathbb{E}\sup\limits_{t\in\left[  0,T\right]  }|X_{t}^{n}|^{2}+\mathbb{E}%
\sup\limits_{t\in\left[  0,T\right]  }|K_{t}^{n}|^{2}+\mathbb{E}%
\updownarrow\!\left\updownarrow K^{n}\right\updownarrow \!\updownarrow_{\ast
T}\leq C_{0}\left[  1+\mathbb{E}|\xi|^{2}+\mathbb{E}|M_{T}^{n}|^{2}\right]
\]
and, if $(\tilde{X}^{n}\left(  \omega\right)  ,\tilde{K}^{n}\left(
\omega\right)  )=\mathcal{GSP}(A;\tilde{\xi}(\omega),\tilde{M}^{n}(\omega))$,
then%
\[
\mathbb{E}\sup\limits_{t\in\left[  0,T\right]  }|X_{t}^{n}-\tilde{X}_{t}%
^{n}|^{2}+\mathbb{E}\sup\limits_{t\in\left[  0,T\right]  }|K_{t}^{n}-\tilde
{K}_{t}^{n}|^{2}\leq C\left(  T\right)  \left[  \mathbb{E}|\xi-\tilde{\xi
}|^{2}+\mathbb{E}|M_{T}^{n}-\tilde{M}_{T}^{n}|^{2}\right]  .
\]
So (replacing $\tilde{M}^{n}$ by $\tilde{M}^{n^{\prime}}$), there exist
$X,K\in L_{ad}^{2}\left(  \Omega;C\left(  \left[  0,T\right]  ;\mathbb{H}%
\right)  \right)  $ such that $X^{n}\rightarrow X$ and $K^{n}\rightarrow K$ in
$L_{ad}^{2}\left(  \Omega;C\left(  \left[  0,T\right]  ;\mathbb{H}\right)
\right)  $ as $n\rightarrow\infty$. The inequalities (\ref{1 SDE-AN}) and
(\ref{2 SDE-AN}) are immediate consequences and, as a by-product, $\left(
X,K\right)  $ is a solution of Eq.(\ref{SDE-AN}).\newline For more details, we
invite the interested reader to consult R\u{a}\c{s}canu \cite{R}.\hfill
\end{proof}

\subsection{Backward stochastic $\mathcal{A}-$representation}

Let $\left(  \Omega,\mathcal{F},\mathbb{P},\{\mathcal{F}_{t}\}_{t\geq
0}\right)  $ be a stochastic basis, where $\{\mathcal{F}_{t}\}_{t\geq0}$ is
the standard filtration associated to a $\mathbb{H}_{0}$-Wiener process
$\{B_{t}\}_{t\geq0}.$

By the representation theorem, for $\xi\in L^{2}\left(  \Omega,\mathcal{F}%
_{T},\mathbb{P};\mathbb{H}\right)  $ there exists a unique $Z\in
\Lambda_{\mathbb{H}\times\mathbb{H}_{0}}^{2}\left(  0,T\right)  $ such that%
\[
\xi=\mathbb{E}\xi+\int_{0}^{T}Z_{s}dB_{s}%
\]
and, for each $\left(  \xi,H\right)  \in L^{2}\left(  \Omega,\mathcal{F}%
_{T},\mathbb{P};\mathbb{H}\right)  \times\Lambda_{\mathbb{H}}^{2}\left(
0,T\right)  $, there exists a unique pair%
\[
\left(  Y,Z\right)  \in S_{\mathbb{H}}^{2}\left[  0,T\right]  \times
\Lambda_{\mathbb{H}\times\mathbb{H}_{0}}^{2}\left(  0,T\right)
\]
such that%
\[
Y_{t}+%
{\displaystyle\int_{t}^{T}}
H_{s}ds=\xi-%
{\displaystyle\int_{t}^{T}}
Z_{s}dB_{s}%
\]
and the mapping $\left(  \xi,H\right)  \mapsto\left(  Y,Z\right)
:L^{2}\left(  \Omega,\mathcal{F}_{T},\mathbb{P};\mathbb{H}\right)
\times\Lambda_{\mathbb{H}}^{2}\left(  0,T\right)  \rightarrow S_{\mathbb{H}%
}^{2}\left[  0,T\right]  \times\Lambda_{\mathbb{H}\times\mathbb{H}_{0}}%
^{2}\left(  0,T\right)  $ is linear and continuous. $\left(  Y,Z\right)  $ is
defined as%
\[
Y_{t}=\mathbb{E}\left(  \left.  \xi-\int_{t}^{T}H_{s}ds\right\vert
\mathcal{F}_{t}\right)  \quad\text{and}\quad\xi-\int_{0}^{T}H_{s}%
ds=\mathbb{E}\left(  \xi-\int_{0}^{T}H_{s}ds\right)  +\int_{0}^{T}Z_{s}%
dB_{s}\text{.}%
\]
Denote%
\[
Y_{t}=C_{t}\left(  \xi,H\right)  \quad\text{and}\quad Z_{t}=D_{t}\left(
\xi,H\right)  .
\]
Remark that, by the Energy Equality, we have%
\begin{equation}
\mathbb{E}\left\vert Y_{t}\right\vert ^{2}+\mathbb{E}%
{\displaystyle\int_{t}^{T}}
\Vert Z_{s}\Vert_{HS}^{2}ds=\mathbb{E}\left\vert \xi\right\vert ^{2}%
+2\mathbb{E}%
{\displaystyle\int_{t}^{T}}
\left\langle Y_{s},F_{s}\right\rangle ds. \label{bsee}%
\end{equation}

If $A:\mathbb{H}\rightrightarrows\mathbb{H}$ is a maximal monotone operator
then the realization of $A$ on $\Lambda_{\mathbb{H}}^{2}\left(  0,T\right)  $
is the maximal monotone operator $\mathcal{A}:\Lambda_{\mathbb{H}}^{2}\left(
0,T\right)  \rightrightarrows\Lambda_{\mathbb{H}}^{2}\left(  0,T\right)  $
defined by $H\in\mathcal{A}\left(  Y\right)  $ iff $H_{t}\left(
\omega\right)  \in A\left(  Y_{t}\left(  \omega\right)  \right)  ,$
$d\mathbb{P}\otimes dt$-$a.e.$ $\left(  \omega,t\right)  \in\Omega
\times\left]  0,T\right[  .$ The inner product in $\Lambda_{\mathbb{H}}%
^{2}\left(  0,T\right)  $ is given by $\left\langle \!\left\langle
U,V\right\rangle \!\right\rangle =\mathbb{E}\int_{0}^{T}\left\langle
U_{t},V_{t}\right\rangle dt.$\medskip

Consider the backward stochastic differential equation%
\begin{equation}
\left\{
\begin{array}
[c]{l}%
-dY_{t}+A\left(  Y_{t}\right)  dt\ni-Z_{t}dB_{t},\ t\in\left[  0,T\right]
,\medskip\\
Y_{T}=\xi,
\end{array}
\right.  \label{A-bsr}%
\end{equation}
where%
\[
\left\{
\begin{array}
[c]{rl}%
\left(  i\right)  & A:\mathbb{H}\rightrightarrows\mathbb{H}\text{ is a maximal
monotone operator and}\medskip\\
\left(  ii\right)  & \xi\in L^{2}(\Omega,\mathcal{F}_{T},\mathbb{P}%
;\overline{Dom(A)}).
\end{array}
\right.
\]

\begin{definition}
$Y\in S_{\mathbb{H}}^{2}\left[  0,T\right]  $ is a solution of Eq.(\ref{A-bsr}%
) if there exist $H\in\Lambda_{\mathbb{H}}^{2}\left(  0,T\right)  $ and
$Z\in\Lambda_{\mathbb{H\times H}_{0}}^{2}\left(  0,T\right)  $ such that%
\[
Y_{t}+%
{\displaystyle\int_{t}^{T}}
H_{s}ds=\xi-%
{\displaystyle\int_{t}^{T}}
Z_{s}dB_{s}%
\]
and $H\in\mathcal{A}\left(  Y\right)  $ (that is, $H_{t}\left(  \omega\right)
\in A\left(  Y_{t}\left(  \omega\right)  \right)  ,$ $d\mathbb{P}\otimes
dt$-$a.e.$ $\left(  \omega,t\right)  \in\Omega\times\left]  0,T\right[  $).
\end{definition}

\noindent Let $R>0$ and the ball $\mathbb{F}_{R}=\left\{  \eta\in L^{2}\left(
\Omega,\mathcal{F}_{T},\mathbb{P};\mathbb{H}\right)  :\mathbb{E}\left\vert
\eta\right\vert ^{2}\leq R\right\}  $.

\noindent For $\left(  U,U^{\ast}\right)  \in\mathcal{A}$ and $\zeta
\in\mathbb{F}_{R}\mathbb{\ }$define%
\[
J_{\left(  \zeta,U,U^{\ast}\right)  }:L^{2}\left(  \Omega,\mathcal{F}%
_{T},\mathbb{P};\mathbb{H}\right)  \mathbb{\times}\Lambda_{\mathbb{H}}%
^{2}\left(  0,T\right)  \times\Lambda_{\mathbb{H}}^{2}\left(  0,T\right)
\rightarrow\mathbb{R}%
\]
by%
\[%
\begin{array}
[c]{rr}%
J_{\left(  \zeta,U,U^{\ast}\right)  }\left(  \eta,Y,H\right)  & =\dfrac{1}%
{2}\mathbb{E}\left\vert \eta-\xi\right\vert ^{2}+\mathbb{E}%
{\displaystyle\int_{0}^{T}}
\left[  \left\langle U_{t},H_{t}\right\rangle +\left\langle Y_{t},U_{t}^{\ast
}\right\rangle -\left\langle U_{t},U_{t}^{\ast}\right\rangle -\left\langle
Y_{t},H_{t}\right\rangle \right]  dt\\
& +\dfrac{1}{2}\left[  \mathbb{E}\left\vert \zeta-\eta\right\vert
^{2}-\mathbb{E}\left\vert \zeta-\xi\right\vert ^{2}\right]
\end{array}
\]
and $\hat{J}:L^{2}\left(  \Omega,\mathcal{F}_{T},\mathbb{P};\mathbb{H}\right)
\mathbb{\times}\Lambda_{\mathbb{H}}^{2}\left(  0,T\right)  \times
\Lambda_{\mathbb{H}}^{2}\left(  0,T\right)  \rightarrow]-\infty,+\infty],$%
\begin{align}
\hat{J}\left(  \eta,Y,H\right)   &  =\sup\left\{  J_{\left(  \zeta,U,U^{\ast
}\right)  }\left(  \eta,Y,H\right)  :\left(  U,U^{\ast}\right)  \in
\mathcal{A},\;\zeta\in\mathbb{F}_{R}\right\}  \medskip
\label{back A representation functional}\\
&  =\dfrac{1}{2}\mathbb{E}\left\vert \eta-\xi\right\vert ^{2}+\mathcal{H}%
\left(  Y,H\right)  -\left\langle \!\left\langle Y,H\right\rangle
\!\right\rangle +\dfrac{1}{2}\sup_{\zeta\in\mathbb{F}_{R}}\left[
\mathbb{E}\left\vert \zeta-\eta\right\vert ^{2}-\mathbb{E}\left\vert \zeta
-\xi\right\vert ^{2}\right]  ,\nonumber
\end{align}
where $\mathcal{H}:\Lambda_{\mathbb{H}}^{2}\left(  0,T\right)  \times
\Lambda_{\mathbb{H}}^{2}\left(  0,T\right)  \rightarrow]-\infty,+\infty]$ is
the Fitzpatrick function associated to the maximal monotone operator
$\mathcal{A}$.

\begin{remark}
$\hat{J}:L^{2}\left(  \Omega,\mathcal{F}_{T},\mathbb{P};\mathbb{H}\right)
\mathbb{\times}\Lambda_{\mathbb{H}}^{2}\left(  0,T\right)  \times
\Lambda_{\mathbb{H}}^{2}\left(  0,T\right)  \rightarrow]-\infty,+\infty]$ is a
l.s.c. function as the supremum of the continuous functions $J_{\left(
\zeta,U,U^{\ast}\right)  }\left(  \eta,Y,H\right)  $.
\end{remark}

\noindent If $\xi\in\mathbb{F}_{R}$ then%
\[
2R^{2}+2\mathbb{E}\left\vert \eta\right\vert ^{2}\geq\sup_{\zeta\in
\mathbb{F}_{R}}\left(  \mathbb{E}\left\vert \zeta-\eta\right\vert
^{2}-\mathbb{E}\left\vert \zeta-\xi\right\vert ^{2}\right)  \geq
\mathbb{E}\left\vert \eta-\xi\right\vert ^{2}%
\]
and clearly follows

\begin{proposition}
\ Let $R>0$ and $\xi\in\mathbb{F}_{R}$. $\hat{J}$ has the following properties:

\begin{itemize}
\item[$\left(  a\right)  $] $\;\hat{J}\left(  \eta,Y,H\right)  \geq
\mathcal{H}\left(  Y,H\right)  -\left\langle \!\left\langle Y,H\right\rangle
\!\right\rangle \geq0$, for all $\left(  \eta,Y,H\right)  \in L^{2}\left(
\Omega,\mathcal{F}_{T},\mathbb{P};\mathbb{H}\right)  \mathbb{\times}%
\Lambda_{\mathbb{H}}^{2}\left(  0,T\right)  \times\Lambda_{\mathbb{H}}%
^{2}\left(  0,T\right)  .$

\item[$\left(  b\right)  $] $\;$Let $(\hat{\eta},\hat{Y},\hat{H})\in
\mathbb{F}_{R}\mathbb{\times}\Lambda_{\mathbb{H}}^{2}\left(  0,T\right)
\times\Lambda_{\mathbb{H}}^{2}\left(  0,T\right)  .$ Then $\hat{J}(\hat{\eta
},\hat{Y},\hat{H})=0$ iff $\hat{\eta}=\xi,$ $\hat{H}\in\mathcal{A(}\hat{Y}).$

\item[$\left(  c\right)  $] $\;$The restriction of $\hat{J}$ to the closed
convex set%
\[
\mathbb{K}=\left\{  \left(  \eta,Y,H\right)  \in\mathbb{F}_{R}\mathbb{\times
}\Lambda_{\mathbb{H}}^{2}\left(  0,T\right)  \times\Lambda_{\mathbb{H}}%
^{2}\left(  0,T\right)  :\;Y_{t}=C_{t}\left(  \eta,H\right)  ,\ \forall
t\in\left[  0,T\right]  \right\}
\]
is a convex lower semicontinuous function and for $(\hat{\eta},\hat{Y},\hat
{H})\in\mathbb{K}$ the following assertions are equivalent:

\begin{itemize}
\item[$\left(  c_{1}\right)  $] $\;\inf\limits_{\left(  \eta,Y,H\right)
\in\mathbb{F}_{R}\mathbb{\times}\Lambda_{\mathbb{H}}^{2}\left(  0,T\right)
\times\Lambda_{\mathbb{H}}^{2}\left(  0,T\right)  }\hat{J}\left(
\eta,Y,H\right)  =\hat{J}(\hat{\eta},\hat{Y},\hat{H})=0.$

\item[$\left(  c_{2}\right)  $] $\;\hat{\eta}=\xi$ and $(\hat{Y},\hat{H}%
,\hat{Z}),$ with $\hat{Z}_{s}=D_{s}(\xi,\hat{H}),$ is the solution of the BSDE
(\ref{A-bsr}).
\end{itemize}
\end{itemize}
\end{proposition}

\begin{proof}
(Sketch) Since the points $\left(  a\right)  $ and $\left(  b\right)  $ are
obvious, we focus on $\left(  c\right)  $. The convexity of $\hat{J}$ on
$\mathbb{K}$ is obtained as follows. By Energy Equality we have%
\begin{gather*}
\frac{1}{2}\left\vert C_{0}\left(  \eta,H\right)  -C_{0}\left(  \zeta
,0\right)  \right\vert ^{2}+\mathbb{E}%
{\displaystyle\int_{0}^{T}}
\left\langle Y_{s}-C_{s}\left(  \zeta,0\right)  ,H_{s}\right\rangle
ds+\frac{1}{2}\mathbb{E}%
{\displaystyle\int_{0}^{T}}
\left\vert D_{s}\left(  \eta,H\right)  -D_{s}\left(  \zeta,0\right)
\right\vert ^{2}ds\\
\quad\quad\quad\quad\quad\quad\quad\quad\quad\quad\quad\quad\quad\quad
\quad\quad\quad\quad\quad\quad\quad\quad\quad\quad=\frac{1}{2}\mathbb{E}%
\left\vert \eta-\zeta\right\vert ^{2}.
\end{gather*}
Then%
\begin{align*}
&  J_{\left(  \zeta,U,\tilde{U}\right)  }\left(  \eta,Y,H\right) \\
&  =\dfrac{1}{2}\mathbb{E}\left\vert \eta-\xi\right\vert ^{2}+\mathbb{E}%
{\displaystyle\int_{0}^{T}}
\left[  \left\langle U_{t},H_{t}\right\rangle +\left\langle Y_{t},U_{t}^{\ast
}\right\rangle -\left\langle U_{t},U_{t}^{\ast}\right\rangle -\left\langle
Y_{t},H_{t}\right\rangle \right]  dt+\dfrac{1}{2}\left[  \mathbb{E}\left\vert
\zeta-\eta\right\vert ^{2}-\mathbb{E}\left\vert \zeta-\xi\right\vert
^{2}\right] \\
&  =\dfrac{1}{2}\mathbb{E}\left\vert \eta-\xi\right\vert ^{2}+\left[
\left\langle \!\left\langle U,H\right\rangle \!\right\rangle +\left\langle
\!\left\langle Y,U^{\ast}\right\rangle \!\right\rangle -\left\langle
\!\left\langle U,U^{\ast}\right\rangle \!\right\rangle \right]  +\frac{1}%
{2}\left\vert C_{0}\left(  \eta,H\right)  -C_{0}\left(  \zeta,0\right)
\right\vert ^{2}\\
&  \quad+\left\langle \!\left\langle C\left(  \zeta,0\right)  ,H\right\rangle
\!\right\rangle +\frac{1}{2}\left\Vert D\left(  \eta,H\right)  -D\left(
\zeta,0\right)  \right\Vert ^{2}-\mathbb{E}\left\vert \zeta-\xi\right\vert
^{2}%
\end{align*}
Hence%
\[%
\begin{array}
[c]{lllll}%
\left(  \eta,Y,H\right)  & \longmapsto & \hat{J}\left(  \eta,Y,H\right)  & = &
\dfrac{1}{2}\mathbb{E}\left\vert \eta-\xi\right\vert ^{2}+\mathcal{H}\left(
Y,H\right)  +\sup\limits_{\zeta}\left\{  \dfrac{1}{2}\left\vert C_{0}\left(
\eta,H\right)  -C_{0}\left(  \zeta,0\right)  \right\vert ^{2}\right.
\medskip\\
&  &  &  & \left.  +\left\langle \!\left\langle C\left(  \zeta,0\right)
,H\right\rangle \!\right\rangle +\dfrac{1}{2}\left\Vert D\left(
\eta,H\right)  -D\left(  \zeta,0\right)  \right\Vert ^{2}-\mathbb{E}\left\vert
\zeta-\xi\right\vert ^{2}\right\}
\end{array}
\]
is, clearly, a convex lower semicontinuous function. Then, the equivalence
between $\left(  c_{1}\right)  $ and $\left(  c_{2}\right)  $ easily
follows.\hfill
\end{proof}

\bigskip

Proving the existence of a solution for the backward stochastic differential
equation (\ref{A-bsr}) is therefore equivalent to solving a problem on convex
analysis. More precisely, it is sufficient to show that the functional defined
by the formula (\ref{back A representation functional}) attains a minimum and
its value in that point is zero. Unfortunately, this is still an open problem,
but we estimate that the perspective and the tools introduced along this paper
will lead us to the desired result.

\section{Fitzpatrick type method for SVI and BSVI}

In the following sections we will consider the finite dimensional case
$\mathbb{H=R}^{d}$ and $\mathbb{H}_{0}=\mathbb{R}^{k}$. Let $\{B_{t},t\geq0\}$
be a $k$-dimensional Brownian motion with respect to a given complete
stochastic basis $\left(  \Omega,\mathcal{F},P,\{\mathcal{F}_{t}\}_{t\geq
0}\right)  $.

\subsection{Stochastic variational inequality}

\subsubsection{Known results}

Let%
\[
F:\Omega\times\left[  0,+\infty\right[  \times\mathbb{R}^{d}\rightarrow
\mathbb{R}^{d},\;\;G:\Omega\times\left[  0,+\infty\right[  \times
\mathbb{R}^{d}\rightarrow\mathbb{R}^{d\times k}.
\]

Consider the stochastic variational inequality (for short SVI)%
\begin{equation}
\left\{
\begin{array}
[c]{l}%
dX_{t}+\partial\varphi(X_{t})(dt)\ni F(t,X_{t})dt+G(t,X_{t})dB_{t},\quad
t\geq0,\smallskip\smallskip\\
X_{0}=\xi,
\end{array}
\right.  \label{MMeq}%
\end{equation}
where will assume%
\begin{equation}%
(%
\mathbf{H}_{0}%
):%
\quad\xi\in L^{0}(\Omega,\mathcal{F}_{0},P;\overline{Dom(\varphi)})
\label{MM-h0}%
\end{equation}
and%
\begin{equation}
(\mathbf{H}_{\varphi}):\quad\left\{
\begin{array}
[c]{rl}%
(i) & \varphi:\mathbb{R}^{d}\rightarrow]-\infty,+\infty]\text{ is a convex
l.s.c. function,}\smallskip\smallskip\\
(ii) & int(Dom(\varphi))\neq\emptyset.
\end{array}
\right.  \label{MM-ha}%
\end{equation}

\begin{definition}
\label{defsub1}A pair $\left(  X,K\right)  \in S_{d}^{0}\times S_{d}^{0}$~,
$K_{0}=0$, is a solution of the stochastic variational inequality (\ref{MMeq})
if the following conditions are satisfied,\
$\mathbb{P}$-%
$%
a.s.%
$:%
\begin{equation}
\left\{
\begin{array}
[c]{ll}%
\left(
d%
_{1}\right)
\;%
& X_{t}\in Dom(\varphi),%
\ a.e.\
t>0\text{ and }\varphi(X)\in L_{loc}^{1}\left(  0,\infty\right)
\text{,}\smallskip\smallskip\\
\left(
d%
_{2}\right)
\;%
& \left\updownarrow K\right\updownarrow _{T}<\infty,\;\forall T>0\text{,}%
\smallskip\smallskip\\
\left(
d%
_{3}\right)
\;%
& X_{t}+K_{t}=\xi+%
{\displaystyle\int_{0}^{t}}
F(s,X_{s})ds+%
{\displaystyle\int_{0}^{t}}
G(s,X_{s})dB_{s},\;\forall t\geq0,\smallskip\smallskip\\
\left(
d%
_{4}\right)
\;%
&
{\displaystyle\int_{s}^{t}}
\left\langle y(r)-X_{r},dK_{r}\right\rangle +%
{\displaystyle\int_{s}^{t}}
\varphi(X_{r})dr\leq%
{\displaystyle\int_{s}^{t}}
\varphi(y(r))dr,\smallskip\\
& \quad\quad\quad\forall y:\mathbb{R}_{+}\rightarrow\mathbb{R}^{d}\text{
continuous function and }\forall0\leq s\leq t.\smallskip\smallskip
\end{array}
\right.  \label{sub2}%
\end{equation}

\end{definition}

\begin{notation}
The notation $dK_{t}\in\partial\varphi\left(  X_{t}\right)  \left(  dt\right)
$ will be used to say that $\left(  X,K\right)  $ satisfy $\left(
d_{1}\right)  ,\left(  d_{2}\right)  $ and $\left(  d_{4}\right)  .$ The SDE
(\ref{MMeq}) will be written, also, in the form%
\[
\left\{
\begin{array}
[c]{l}%
X_{t}+K_{t}=\xi+%
{\displaystyle\int_{0}^{t}}
F(s,X_{s})ds+%
{\displaystyle\int_{0}^{t}}
G(s,X_{s})dB_{s},\text{ }\forall t\geq0,\smallskip\smallskip\\
dK_{t}\in\partial\varphi\left(  X_{t}\right)  \left(  dt\right)  .
\end{array}
\right.
\]

\end{notation}

Remark (see Asiminoaei \& R\u{a}\c{s}canu \cite{AR}) that the condition
$\left(  d_{4}\right)  $ from Definition \ref{defsub1} is equivalent to each
of the following conditions, for any fixed $T>0$,%
\[%
\begin{array}
[c]{ll}%
(a_{1}) &
{\displaystyle\int\nolimits_{s}^{t}}
\left\langle z-X_{r},dK_{r}\right\rangle +%
{\displaystyle\int\nolimits_{s}^{t}}
\varphi(X_{r})dr\leq(t-s)\varphi(z),\ \forall z\in\mathbb{R}^{d}%
,\;\forall0\leq s\leq t\leq T,\smallskip\\
(a_{2}) &
{\displaystyle\int\nolimits_{s}^{t}}
\left\langle X_{r}-z,dK_{r}-z^{\ast}dr\right\rangle \geq0,\,\,\forall\,\left(
z,z^{\ast}\right)  \in\partial\varphi,\ \forall0\leq s\leq t\leq
T,\smallskip\\
(a_{3}) &
{\displaystyle\int\nolimits_{0}^{T}}
\left\langle y(r)-X_{r},dK_{r}\right\rangle +%
{\displaystyle\int\nolimits_{0}^{T}}
\varphi(X_{r})dr\leq%
{\displaystyle\int\nolimits_{0}^{T}}
\varphi(y(r))dr,\ \forall y\in C([0,T],\mathbb{R}^{d}).\smallskip
\end{array}
\]

\noindent Hence, the condition $\left(
d%
_{4}\right)  $ means that $\left(  X_{\cdot}\left(  \omega\right)  ,K_{\cdot
}\left(  \omega\right)  \right)  \in\partial\tilde{\varphi},$ $\mathbb{P}%
$-a.s., where $\tilde{\varphi}$ is the realization of $\varphi$ on $C\left(
\left[  0,T\right]  ;\mathbb{R}^{d}\right)  $, that is $\tilde{\varphi
}:C([0,T];\mathbb{R}^{d})\rightarrow]-\infty,+\infty],$
\begin{equation}
\tilde{\varphi}(x)=\left\{
\begin{array}
[c]{l}%
{\displaystyle\int\nolimits_{0}^{T}}
\varphi(x(t))dt,\,\,\text{\textit{if }}\varphi\left(  x\right)  \in
L^{1}(0,T),\smallskip\smallskip\\
+\infty,\;\text{\textit{otherwise.}}%
\end{array}
\right.  \label{Ba2}%
\end{equation}

\begin{notation}
We introduce the notation:%
\[
F_{R}^{\#}\left(  t\right)  :=\operatorname*{ess}\sup\left\{  \left\vert
F(t,x)\right\vert :\left\vert x\right\vert \leq R\right\}  .
\]

\end{notation}

We recall the basic assumptions on $F$ and $G$ under which we will study the
multivalued stochastic equation (\ref{MMeq}):\smallskip\smallskip
\newline$\circ\quad$the functions $F\left(  \cdot,\cdot,x\right)
:\Omega\times\left[  0,+\infty\right[  \rightarrow\mathbb{R}^{d}$ and
$G\left(  \cdot,\cdot,x\right)  :\Omega\times\left[  0,+\infty\right[
\rightarrow\mathbb{R}^{d\times k}$ are progressively measurable stochastic
processes for every $x\in\mathbb{R}^{d},$\smallskip\smallskip\newline%
$\circ\quad$there exist $\mu\in L_{loc}^{1}\left(  0,\infty\right)  $ and
$\ell\in L_{loc}^{2}\left(  0,\infty;\mathbb{R}_{+}\right)  $, such that
$d\mathbb{P}\otimes dt$-$a.e.$:%
\begin{equation}
\left(  \mathbf{H}_{F}\right)  :\left\{
\begin{array}
[c]{ll}
& \text{Continuity:}\smallskip\\
\left(  \mathbf{C}_{F}\right)  :\quad & x\mapsto F\left(  t,x\right)
:\mathbb{R}^{d}\rightarrow\mathbb{R}^{d}\text{ is continuous,}\smallskip
\smallskip\\
& \text{Monotonicity condition:}\smallskip\\
\left(  \mathbf{M}_{F}\right)  :\quad & \left\langle
x-y,F(t,x)-F(t,y)\right\rangle \,\leq\mu\left(  t\right)  |x-y|^{2},\text{
}\forall x,y\in\mathbb{R}^{d}\text{,}\smallskip\smallskip\\
& \text{Boundedness condition:}\smallskip\\
\left(  \mathbf{B}_{F}\right)  :\quad &
{\displaystyle\int_{0}^{T}}
F_{R}^{\#}\left(  s\right)  ds<\infty\text{, for all }R,T\geq0.
\end{array}
\right.  \label{MM-hf}%
\end{equation}
and
\begin{equation}
(\mathbf{H}_{G})~:\left\{
\begin{array}
[c]{ll}
& \text{Lipschitz condition:}\smallskip\\
\left(  \mathbf{L}_{G}\right)  :\quad & |G(t,x)-G(t,y)|\leq\ell\left(
t\right)  |x-y|,\text{ }\forall x,y\in\mathbb{R}^{d}\text{,}\smallskip
\smallskip\\
& \text{Boundedness condition:}\smallskip\\
\left(  \mathbf{B}_{g}\right)  :\quad &
{\displaystyle\int_{0}^{T}}
|G(t,0)|^{2}dt<\infty.
\end{array}
\right.  \label{MM-hg}%
\end{equation}
Clearly $\left(  \mathbf{H}_{F}\right)  $ and $(\mathbf{H}_{G})$ yield
$F(\cdot,\cdot,X_{\cdot})\in L_{loc}^{1}\left(  \mathbb{R}_{+};\mathbb{R}%
^{d}\right)  \;$and $G(\cdot,\cdot,X_{\cdot})\in\Lambda_{d\times k}^{0}$ for
all $X\in S_{d}^{0}$.

\begin{theorem}
\label{MM-t1-exist}If the assumptions (\ref{MM-h0}), (\ref{MM-ha}),
(\ref{MM-hf}) and (\ref{MM-hg}) are satisfied, then the SDE (\ref{MMeq}) has a
unique solution $\left(  X,K\right)  \in S_{d}^{0}\times S_{d}^{0}$ (in the
sense of Definition \ref{defsub1}). Moreover, if there exist $p\geq2$ and
$u_{0}\in int\left(  Dom\left(  \varphi\right)  \right)  $ such that, for all
$T\geq0$,%
\begin{equation}
\mathbb{E}\left\vert \xi\right\vert ^{p}+\mathbb{E}\left(
{\displaystyle\int_{0}^{T}}
\left\vert F\left(  t,u_{0}\right)  \right\vert dt\right)  ^{p}+\mathbb{E}%
\left(
{\displaystyle\int_{0}^{T}}
|G(t,u_{0})|^{2}dt\right)  ^{p/2}<+\infty\text{,} \label{MM-7-a}%
\end{equation}
then%
\[
\mathbb{E}(\left\Vert X\right\Vert _{T}^{p}+\left\Vert K\right\Vert _{T}%
^{p/2}+\left\updownarrow K\right\updownarrow _{T}^{p/2})+\mathbb{E}\left(
{\displaystyle\int_{0}^{T}}
\left\vert \varphi\left(  X_{r}\right)  \right\vert dr\right)  ^{p/2}<\infty.
\]

\end{theorem}

(For the proof see Pardoux \& R\u{a}\c{s}canu \cite{PRb}, Theorem
4.14.)\bigskip

\subsubsection{Fitzpatrick approach}

In this subsection, assumptions $\left(  \mathbf{H}_{F}\right)  $ and $\left(
\mathbf{H}_{G}\right)  $ are replaced by\smallskip\smallskip\newline$\left(
i\right)  \quad$the functions $F\left(  \cdot,\cdot,x\right)  :\Omega
\times\left[  0,+\infty\right[  \rightarrow\mathbb{R}^{d}$ and $G\left(
\cdot,\cdot,x\right)  :\Omega\times\left[  0,+\infty\right[  \rightarrow
\mathbb{R}^{d\times k}$ are progressively measurable stochastic processes for
every $x\in\mathbb{R}^{d}$ and, $d\mathbb{P}\otimes dt$-$a.e.$,\smallskip
\smallskip\newline$\left(  ii\right)  \quad x\mapsto F\left(  t,x\right)
:\mathbb{R}^{d}\rightarrow\mathbb{R}^{d}$ and $x\mapsto G\left(  t,x\right)
:$ $\mathbb{R}^{d}\rightarrow\mathbb{R}^{d\times k}$ are $\,$%
continuous,\smallskip\smallskip\newline$\left(  iii\right)  \quad$for all
$x,y\in\mathbb{R}^{d}$%
\begin{equation}
2\left\langle x-y,F(t,x)-F(t,y)\right\rangle +|G(t,x)-G(t,y)|^{2}\leq0\text{
and} \label{svi-a1}%
\end{equation}
$\left(  iv\right)  $ there exists $b>0$ such that, for all $x\in
\mathbb{R}^{d}$,%
\begin{equation}
\left\vert F(t,x)\right\vert +\left\vert G(t,x)\right\vert \leq b\left(
1+\left\vert x\right\vert \right)  . \label{svi-a2}%
\end{equation}

\begin{remark}
If $\mu\left(  t\right)  +\frac{1}{2}\ell^{2}\left(  t\right)  \leq0$,$\ $for
every $t\geq0$, then the assumptions (\ref{MM-hf}-$\mathbf{M}_{F}$) and
(\ref{MM-hg}-$\mathbf{L}_{G}$) implies that (\ref{svi-a1}) holds.
\end{remark}

Denote%
\[
\mathbb{S}_{BV}\left[  0,T\right]  =\left\{  K\in S_{d}^{0}\left[  0,T\right]
:K_{0}=0,\;\mathbb{E}\left\updownarrow K\right\updownarrow _{T}^{2}%
<\infty\right\}  ,
\]
with the $w^{\ast}$-topology, that means $K^{n}\rightarrow K$ if
$\lim\limits_{n\rightarrow\infty}\mathbb{E}%
{\textstyle\int_{0}^{T}}
\left\langle X_{t},dK_{t}^{n}\right\rangle =\mathbb{E}%
{\textstyle\int_{0}^{T}}
\left\langle X_{t},dK_{t}\right\rangle $, for all $X\in L^{2}(\Omega
;C([0,T];\mathbb{R}^{d}))$.

Let $\Phi:$\ $S_{d}^{2}\left[  0,T\right]  \rightarrow]-\infty,+\infty]$
defined by
\begin{equation}%
\begin{tabular}
[c]{lll}%
$\Phi(X)$ & $=$ & $\left\{
\begin{array}
[c]{l}%
\mathbb{E}%
{\displaystyle\int\nolimits_{0}^{T}}
\varphi(X_{t})dt,\;\text{\textit{if }}\varphi\left(  X\right)  \in
L^{1}(\Omega\times]0,T[),\smallskip\vspace{0.04in}\\
+\infty,\;\text{\textit{otherwise.}}%
\end{array}
\right.  $%
\end{tabular}
\ \ \ \ \ \ \ \ \ \label{svi-FI}%
\end{equation}
Since\ $\varphi:\mathbb{R}^{d}\rightarrow]-\infty,+\infty]$ is a proper convex
l.s.c. function then $\Phi$ is also a proper convex l.s.c.~function.

Let%
\[
\mathbb{S}:=L^{2}(\Omega,\mathcal{F}_{0},\mathbb{P};\overline{Dom(\varphi
)})\times S_{d}^{2}\left[  0,T\right]  \times\mathbb{S}_{BV}\left[
0,T\right]  \times\Lambda_{d\times k}^{2}\left(  0,T\right)
\]
and, for each $U\in Dom\left(  \Phi\right)  =\left\{  X\in S_{d}^{2}\left[
0,T\right]  :\Phi(X)<\infty\right\}  $, we consider the mapping\newline%
$J_{U}:\mathbb{S}\rightarrow]-\infty,+\infty]$, defined by%
\begin{equation}%
\begin{array}
[c]{rr}%
J_{U}\left(  \eta,X,L,g\right)  & =\dfrac{1}{2}\mathbb{E}\left\vert \eta
-\xi\right\vert ^{2}+\mathbb{E}%
{\displaystyle\int_{0}^{T}}
\left[  \left\langle U_{s}-X_{s},F\left(  s,U_{s}\right)  \right\rangle
\mathbb{+}\dfrac{1}{2}\left\vert g_{s}-G\left(  s,U_{s}\right)  \right\vert
^{2}\right]  ds\medskip\\
& +\mathbb{E}%
{\displaystyle\int_{0}^{T}}
\left\langle U_{s}-X_{s},dL_{s}\right\rangle +\Phi\left(  X\right)
-\Phi\left(  U\right)
\end{array}
\label{functsvi}%
\end{equation}

and $\hat{J}:\mathbb{S}\longrightarrow]-\infty,+\infty]$%
\[
\hat{J}\left(  \eta,X,L,g\right)  :=\underset{U\in Dom\left(  \Phi\right)
}{\sup}J_{U}\left(  \eta,X,L,g\right)  .
\]

\begin{remark}
$\hat{J}:\mathbb{S}\rightarrow]-\infty,+\infty]$ is a lower semicontinuous
function as supremum of lower semicontinuous functions.
\end{remark}

We now have

\begin{proposition}
$\hat{J}$ has the following properties:

\begin{itemize}
\item[$\left(  a\right)  $] $\;\hat{J}\left(  \eta,X,L,g\right)  \geq0$, for
all $\left(  \eta,X,L,g\right)  \in\mathbb{S}$ and $\hat{J}$ is not
identically $+\infty.$

\item[$\left(  b\right)  $] $\;$Let $(\hat{\eta},\hat{X},\hat{L},\hat{g}%
)\in\mathbb{S}.$ Then%
\[
\hat{J}(\hat{\eta},\hat{X},\hat{L},\hat{g})=0\quad\text{iff}\quad\hat{\eta
}=\xi,\;\hat{g}_{\cdot}=G(\cdot,\hat{X}_{\cdot}),\;\hat{L}+%
{\displaystyle\int_{0}^{\cdot}}
F(s,\hat{X}_{s})ds\in\partial\Phi(\hat{X}).
\]

\item[$\left(  c\right)  $] $\;$The restriction of $\hat{J}$ to the closed
convex set%
\[
\mathbb{L}=\left\{  \left(  \eta,X,L,g\right)  \in\mathbb{S}:X_{t}+L_{t}=\eta+%
{\displaystyle\int_{0}^{t}}
g_{s}dB_{s},\ \forall t\in\left[  0,T\right]  \right\}
\]
is a convex l.s.c. function. If $(\hat{\eta},\hat{X},\hat{L},\hat{g}%
)\in\mathbb{L},$ then $\hat{J}(\hat{\eta},\hat{X},\hat{L},\hat{g})=0$ iff%
\[
\hat{\eta}=\xi,\;\hat{g}_{\cdot}=G(\cdot,\hat{X}_{\cdot})\text{ and }(\hat
{X},\hat{L}+%
{\textstyle\int_{0}^{\cdot}}
F(s,\hat{X}_{s})ds)\text{ is a solution of the SVI (\ref{MMeq}).}%
\]

\end{itemize}
\end{proposition}

\begin{proof}
$\left(  a\right)  \;$If $X\notin Dom\left(  \Phi\right)  $ then $\hat
{J}\left(  \eta,X,L,g\right)  =+\infty.$ If $X\in Dom\left(  \Phi\right)  $
then%
\begin{align*}
\hat{J}\left(  \eta,X,L,g\right)   &  =\underset{U\in Dom\left(  \Phi\right)
}{\sup}J_{U}\left(  \eta,X,L,g\right) \\
&  \geq J_{X}\left(  \eta,X,L,g\right) \\
&  =\dfrac{1}{2}\mathbb{E}\left\vert \eta-\xi\right\vert ^{2}+\dfrac{1}%
{2}\mathbb{E}%
{\displaystyle\int_{0}^{T}}
\left\vert g_{s}-G\left(  s,X_{s}\right)  \right\vert ^{2}ds\\
&  \geq0.
\end{align*}
$\hat{J}$ is a proper function since, for $v_{0}\in\partial\varphi\left(
u_{0}\right)  $ and $\eta^{0}=\xi,$ $X_{t}^{0}=u_{0},\;L_{t}^{0}=v_{0}t-%
{\textstyle\int_{0}^{t}}
F\left(  s,u_{0}\right)  ds,$ $g_{s}^{0}=G\left(  s,u_{0}\right)  $, we have
(using the assumption (\ref{svi-a1})) that%
\[
J_{U}\left(  \eta^{0},X^{0},L^{0},g^{0}\right)  \leq0\text{,}\;\text{for all
}U\in Dom\left(  \Phi\right)  \text{.}%
\]

$\left(  b\right)  \;$If $\hat{J}(\hat{\eta},\hat{X},\hat{L},\hat{g})=0$, then
$\hat{X}\in Dom\left(  \Phi\right)  $ and by the calculus from the proof of
$\left(  a\right)  $ we infer $\hat{\eta}=\xi,$ $\hat{g}=G(\cdot,\hat
{X}_{\cdot})$ and%
\[
J_{U}(\hat{\eta},\hat{X},\hat{L},\hat{g})\leq0,\;\text{for all }U\in
Dom\left(  \Phi\right)  .
\]
Hence%
\[
\mathbb{E}%
{\displaystyle\int_{0}^{T}}
\left\langle U_{s}-\hat{X}_{s},F\left(  s,U_{s}\right)  ds+d\hat{L}%
_{s}\right\rangle +\Phi(\hat{X})\leq\Phi\left(  U\right)  \text{,}\;\text{for
all }U\in Dom\left(  \Phi\right)  .
\]
Let $V\in Dom\left(  \Phi\right)  $ and $\lambda\in]0,1[$ be arbitrary. Since
$Dom\left(  \Phi\right)  $ is a convex set, we can replace $U$ by $\left(
1-\lambda\right)  \hat{X}+\lambda V.$ It follows%
\begin{align*}
&  \lambda\mathbb{E}%
{\displaystyle\int_{0}^{T}}
\left\langle V_{s}-\hat{X}_{s},F(s,\hat{X}_{s}+\lambda(V_{s}-\hat{X}%
_{s}))ds+d\hat{L}_{s}\right\rangle +\Phi(\hat{X})\\
&  \leq\Phi(\left(  1-\lambda\right)  \hat{X}+\lambda V)\leq\left(
1-\lambda\right)  \Phi(\hat{X})+\lambda\Phi\left(  V\right)  \text{,}%
\end{align*}
which is equivalent to%
\[
\mathbb{E}%
{\displaystyle\int_{0}^{T}}
\left\langle V_{s}-\hat{X}_{s},F(s,\hat{X}_{s}+\lambda(V_{s}-\hat{X}%
_{s}))ds+d\hat{L}_{s}\right\rangle +\Phi(\hat{X})\leq\Phi\left(  V\right)
\text{,}%
\]
for all $V\in Dom\left(  \Phi\right)  .$ By the continuity of $x\mapsto
F\left(  t,x\right)  $ and assumption (\ref{svi-a2}) we can pass to limit
under the last integral, and it follows that $\hat{L}+%
{\textstyle\int_{0}^{\cdot}}
F(s,\hat{X}_{s})ds\in\partial\Phi(\hat{X})$.

Conversely, using (\ref{svi-a1}), we have%
\begin{align*}
&  J_{U}(\xi,\hat{X},\hat{L},G(\cdot,\hat{X}_{.}))\\
\quad &  =\dfrac{1}{2}\mathbb{E}%
{\displaystyle\int_{0}^{T}}
|G(s,\hat{X}_{s})-G\left(  s,U_{s}\right)  |^{2}ds+\mathbb{E}%
{\displaystyle\int_{0}^{T}}
\left\langle U_{s}-\hat{X}_{s},F\left(  s,U_{s}\right)  -F(s,\hat{X}%
_{s})ds\right\rangle \\
&  +\mathbb{E}%
{\displaystyle\int_{0}^{T}}
\left\langle U_{s}-\hat{X}_{s},F(s,\hat{X}_{s})ds+d\hat{L}_{s}\right\rangle
+\Phi(\hat{X})-\Phi\left(  U\right) \\
&  \leq0
\end{align*}
and, consequently, $\hat{J}(\xi,\hat{X},\hat{L},G(\cdot,\hat{X}_{.}))=0$.

$\left(  c\right)  \;$The second part of this point is easy to observe, and,
therefore, $(\hat{X},\hat{L}+%
{\textstyle\int_{0}^{\cdot}}
F(s,\hat{X}_{s})ds)$ is a solution of the SVI (\ref{MMeq}).

It remains to prove the convexity of $\hat{J}$ on $\mathbb{L}.$ By the Energy
Equality we have%
\[
\frac{1}{2}\mathbb{E}\left\vert X_{T}\right\vert ^{2}+\mathbb{E}%
{\displaystyle\int_{0}^{T}}
\left\langle X_{s},dL_{s}\right\rangle =\frac{1}{2}\mathbb{E}\left\vert
\eta\right\vert ^{2}+\frac{1}{2}\mathbb{E}%
{\displaystyle\int_{0}^{T}}
\left\vert g_{s}\right\vert ^{2}ds
\]
and, using it in the formula (\ref{functsvi}), the functional $J_{U}\left(
\eta,X,L,g\right)  $ becomes%
\begin{align*}
J_{U}\left(  \eta,X,L,g\right)   &  =\dfrac{1}{2}\mathbb{E}\left\vert \eta
-\xi\right\vert ^{2}+\mathbb{E}%
{\displaystyle\int_{0}^{T}}
\left\langle U_{s}-X_{s},F\left(  s,U_{s}\right)  ds\right\rangle +\dfrac
{1}{2}\mathbb{E}%
{\displaystyle\int_{0}^{T}}
\left\vert g_{s}-G\left(  s,U_{s}\right)  \right\vert ^{2}ds\\
&  +\left[  \mathbb{E}%
{\displaystyle\int_{0}^{T}}
\left\langle U_{s},dL_{s}\right\rangle -\dfrac{1}{2}\mathbb{E}\left\vert
\eta\right\vert ^{2}-\dfrac{1}{2}\mathbb{E}%
{\displaystyle\int_{0}^{T}}
\left\vert g_{s}\right\vert ^{2}ds+\dfrac{1}{2}\mathbb{E}\left\vert
X_{T}\right\vert ^{2}\right]  +\Phi\left(  X\right)  -\Phi\left(  U\right) \\
&  =-\mathbb{E}\left\langle \eta,\xi\right\rangle +\dfrac{1}{2}\mathbb{E}%
\left\vert \xi\right\vert ^{2}+\mathbb{E}%
{\displaystyle\int_{0}^{T}}
\left\langle U_{s}-X_{s},F\left(  s,U_{s}\right)  ds\right\rangle +\mathbb{E}%
{\displaystyle\int_{0}^{T}}
\left\langle U_{s},dL_{s}\right\rangle \\
&  +\dfrac{1}{2}\mathbb{E}\left\vert X_{T}\right\vert ^{2}+\dfrac{1}%
{2}\mathbb{E}%
{\displaystyle\int_{0}^{T}}
\left\vert G\left(  s,U_{s}\right)  \right\vert ^{2}ds-\mathbb{E}%
{\displaystyle\int_{0}^{T}}
\left\langle g_{s},G\left(  s,U_{s}\right)  \right\rangle ds+\Phi\left(
X\right)  -\Phi\left(  U\right)  .
\end{align*}
It clearly follows that $J_{U}$ is convex and lower semicontinuous for
$\forall U\in Dom\left(  \Phi\right)  .$ Consequently, the mapping $\left(
\eta,X,L,g\right)  \longmapsto\hat{J}\left(  \eta,X,L,g\right)  =\sup_{U\in
Dom\left(  \Phi\right)  }J_{U}\left(  \eta,X,L,g\right)  $ has the same properties.

The proof is now complete.\hfill
\end{proof}

\subsection{Backward stochastic variational inequality}

In this section we suppose that the filtration $\left\{  \mathcal{F}%
_{t}:\,t\geq0\right\}  $ is the natural filtration of the $k$--dimensional
Brownian motion $\{B_{t}\,:\,t\geq0\}$, i.e., for all $t\geq0$,%
\[
\mathcal{F}_{t}=\mathcal{F}_{t}^{B}:=\sigma\left(  \left\{  B_{s}:0\leq s\leq
t\right\}  \right)  \vee\mathcal{N}_{\mathbb{P}}.
\]

\subsubsection{Known results}

Consider the backward stochastic variational inequality (for short BSVI)%
\begin{equation}
\left\{
\begin{array}
[c]{l}%
-dY_{t}+\partial\varphi\left(  Y_{t}\right)  dt\ni F\left(  t,Y_{t}%
,Z_{t}\right)  dt-Z_{t}dB_{t},\;0\leq t<T,\medskip\\
Y_{T}=\xi,
\end{array}
\,\right.  \label{bsvi1}%
\end{equation}
or, equivalently,%
\[
\left\{
\begin{array}
[c]{l}%
Y_{t}+%
{\displaystyle\int_{t}^{T}}
H_{s}ds=\xi+%
{\displaystyle\int_{t}^{T}}
F\left(  s,Y_{s},Z_{s}\right)  ds-\int_{t}^{T}Z_{s}dB_{s},\;t\in\left[
0,T\right]  ,\;\mathbb{\mathbb{P}}\text{-}a.s.,\medskip\\
H_{t}\left(  \omega\right)  \in\partial\varphi\left(  Y_{t}\left(
\omega\right)  \right)  ,\text{ }d\mathbb{P}\otimes dt\text{-}a.e.
\end{array}
\right.
\]
We assume

\begin{itemize}
\item[$\left(  \mathbf{H}_{\xi}\right)  $] :$\quad\xi:\Omega\rightarrow
\mathbb{R}^{d}$\textit{ }is\textit{ a }$F_{T}$-measurable random vector,

\item[$\left(  \mathbf{H}_{\varphi}\right)  $] :$\quad\partial\varphi$ is the
subdifferential of the proper convex l.s.c. function $\varphi:\mathbb{R}%
^{d}\rightarrow]-\infty,+\infty]$,

\item[$\left(  \mathbf{H}_{F}\right)  $] :$\quad F:\Omega\times\left[
0,\infty\right[  \times\mathbb{R}^{d}\times\mathbb{R}^{d\times k}%
\rightarrow\mathbb{R}^{d}$\textit{ }satisfies

\begin{itemize}
\item[$\circ$] $\quad$the function $F\left(  \cdot,\cdot,y,z\right)
:\Omega\times\left[  0,T\right]  \rightarrow\mathbb{R}^{d}$ is a progressively
measurable stochastic process for every $\left(  y,z\right)  \in\mathbb{R}%
^{d}\times\mathbb{R}^{d\times k}$,

\item[$\circ$] $\quad$there exist some deterministic functions\textit{ }%
$\mu\in L^{1}\left(  0,T;\mathbb{R}\right)  $ and\textit{ }$\ell\in
L^{2}\left(  0,T;\mathbb{R}\right)  $, such that,%
\begin{equation}
\left\{
\begin{array}
[c]{l}%
\left(  i\right)  \;\text{for all }y,y^{\prime}\in\mathbb{R}^{d}%
,\;z,z^{\prime}\in\mathbb{R}^{d\times k}\text{,}\;d\mathbb{P}\otimes
dt\text{-}a.e.:\smallskip\\
\;\;\;\;\;\;\;\;\;%
\begin{array}
[c]{ll}
& \text{Continuity:}\smallskip\\
\left(  \mathbf{C}_{y}\right)  : & y\longrightarrow F\left(  t,y,z\right)
:\mathbb{R}^{d}\rightarrow\mathbb{R}^{d}\,\text{is continuous,}\smallskip
\smallskip\\
& \text{Monotonicity condition:}\smallskip\\
\left(  \mathbf{M}_{y}\right)  : & \left\langle y^{\prime}-y,F(t,y^{\prime
},z)-F(t,y,z)\right\rangle \leq\mu\left(  t\right)  |y^{\prime}-y|^{2}%
\text{,}\smallskip\smallskip\\
& \text{Lipschitz condition:}\smallskip\\
\left(  \mathbf{L}_{z}\right)  : & |F(t,y,z^{\prime})-F(t,y,z)|\leq\ell\left(
t\right)  ~|z^{\prime}-z|\text{,}\smallskip\smallskip
\end{array}
\\
\left(  ii\right)  \;\text{Boundedness condition:}\\
\;\;\;\;\;\;\;\;\;%
\begin{array}
[c]{ll}%
\left(  \mathbf{B}_{F}\right)  \quad &
{\displaystyle\int_{0}^{T}}
F_{R}^{\#}\left(  t\right)  dt<\infty,\;\mathbb{\mathbb{P}}\text{-}a.s.\text{,
}\forall R\geq0\text{,}%
\end{array}
\end{array}
\right.  \label{Ch5-GM-ip1}%
\end{equation}

\end{itemize}
\end{itemize}

where%
\[
F_{R}^{\#}\left(  t\right)  =\sup\left\{  \left\vert F(t,y,0)\right\vert
:\left\vert y\right\vert \leq R\right\}  .
\]

\begin{definition}
\label{def1-bsvi} A pair $\left(  Y,Z\right)  \in S_{d}^{0}\left[  0,T\right]
\times\Lambda_{d\times k}^{0}\left(  0,T\right)  $ of stochastic processes is
a solution of the backward stochastic variational inequality (\ref{bsvi1}) if
there exists a progressively measurable stochastic process $H$ such that,
$\mathbb{P}$-$a.s.$,%
\[%
\begin{array}
[c]{ll}%
\left(  a\right)  \;\; &
{\displaystyle\int_{0}^{T}}
\left\vert H_{t}\right\vert dt+%
{\displaystyle\int_{0}^{T}}
\left\vert F(t,Y_{t},Z_{t})\right\vert dt<\infty,\medskip\\
\left(  b\right)  \;\; & \left(  Y_{t}\left(  \omega\right)  ,H_{t}\left(
\omega\right)  \right)  \in\partial\varphi,\text{ }a.e.\text{ }t\in\left[
0,T\right]
\end{array}
\]
and, for all $t\in\left[  0,T\right]  $,%
\begin{equation}
Y_{t}+\int_{t}^{T}H_{s}ds=\eta+\int_{t}^{T}F\left(  s,Y_{s},Z_{s}\right)
ds-\int_{t}^{T}Z_{s}dB_{s}. \label{bsvi2}%
\end{equation}
(Without confusion, \textit{the uniqueness of the stochastic process }$H$ will
permit to call the triplet $\left(  Y,Z,H\right)  $ a solution of
Eq.(\ref{bsvi1}).)
\end{definition}

We introduce now a supplementary assumption\medskip

$\left(  \mathbf{A}\right)  :$ There exist $p\geq2,$ a positive stochastic
process $\beta\in L^{1}\left(  \Omega\times\left]  0,T\right[  \right)  ,$ a
positive function $b\in L^{1}\left(  0,T\right)  $ and a real number
$\kappa\geq0$, such that for all\ $\left(  u,\hat{u}\right)  \in
\partial\varphi\ $and\ $z\in R^{d\times k}$%
\[
\left\langle \hat{u},F\left(  t,u,z\right)  \right\rangle \leq\dfrac{1}%
{2}\left\vert \hat{u}\right\vert ^{2}+\beta_{t}+b\left(  t\right)  \left\vert
u\right\vert ^{p}+\kappa\left\vert z\right\vert ^{2},\text{ }d\mathbb{P}%
\otimes dt\text{-}a.e.
\]

\begin{theorem}
\label{t1-bsvi}Let assumptions $\left(  \mathbf{H}_{\xi}\right)  $, $\left(
\mathbf{H}_{\varphi}\right)  $, $\left(  \mathbf{H}_{F}\right)  $ and $\left(
\mathbf{A}\right)  $ be satisfied. If there exists $u_{0}\in Dom\left(
\partial\varphi\right)  $ such that%
\begin{equation}
\mathbb{E}\left\vert \xi\right\vert ^{p}+\mathbb{E}\left\vert \varphi\left(
\xi\right)  \right\vert +\mathbb{E}\left(
{\displaystyle\int_{0}^{T}}
\left\vert F(s,u_{0},0)\right\vert ds\right)  ^{p}<\infty, \label{bsvi-ip2}%
\end{equation}
then the BSVI (\ref{bsvi1}) has a unique solution $\left(  Y,Z\right)  \in
S_{d}^{p}\left[  0,T\right]  \times\Lambda_{d\times k}^{p}\left(  0,T\right)
.$ Moreover, uniqueness holds in $S_{d}^{1+}\left[  0,T\right]  \times
\Lambda_{d\times k}^{0}\left(  0,T\right)  ,$ where%
\[
S_{d}^{1+}\left[  0,T\right]  :=%
{\displaystyle\bigcup_{p>1}}
S_{d}^{p}\left[  0,T\right]  .
\]

\end{theorem}

(For the proof see Pardoux \& R\u{a}\c{s}canu \cite{PRb}, Theorem
5.13.)\bigskip

\subsubsection{Fitzpatrick approach}

In this subsection the assumptions $\left(  \mathbf{H}_{F}\right)  $ are
replaced by\smallskip\smallskip\newline$\left(  i\right)  \quad$the function
$F\left(  \cdot,\cdot,y,z\right)  :\Omega\times\left[  0,+\infty\right[
\rightarrow\mathbb{R}^{d}$ is a progressively measurable stochastic processes
for every $\left(  y,z\right)  \in\mathbb{R}^{d}\times\mathbb{R}^{d\times k}%
$,\smallskip\smallskip\newline$\left(  ii\right)  \quad\left(  y,z\right)
\mapsto F\left(  t,y,z\right)  :\mathbb{R}^{d}\times\mathbb{R}^{d\times
k}\rightarrow\mathbb{R}^{d}$ is $\,$continuous $d\mathbb{P}\otimes dt$%
-$a.e.$,\smallskip\smallskip\newline$\left(  iii\right)  \quad$for all
$y,y^{\prime}\in\mathbb{R}^{d}$ and $z,z^{\prime}\in\mathbb{R}^{d\times k}$%
\begin{equation}
\left\langle y-y^{\prime},F(t,y,z)-F(t,y^{\prime},z^{\prime})\right\rangle
\leq\frac{1}{2}\left\vert z-z^{\prime}\right\vert ,\;d\mathbb{P}\otimes
dt\text{-}a.e., \label{bsvi-a1}%
\end{equation}
$\left(  iv\right)  $ there exists $b>0$ such that, for all $y\in
\mathbb{R}^{d}$,%
\[
\left\vert F(t,y,z)\right\vert \leq b\left(  1+\left\vert y\right\vert
+\left\vert z\right\vert \right)  ,\;d\mathbb{P}\otimes dt\text{-}a.e.
\]

Remark that, if%
\[
\mu\left(  t\right)  +\frac{1}{2}\ell^{2}\left(  t\right)  \leq0,\;a.e.\text{
}t\geq0,
\]
then the assumptions $\left(  \mathbf{H}_{F}\right)  $ implies $\left(
i\right)  -\left(  iii\right)  .$

Denote by $\Phi:S_{d}^{2}\left[  0,T\right]  \rightarrow]-\infty,+\infty]$ the
proper convex lower semicontinuous function defined by%
\[%
\begin{tabular}
[c]{lll}%
$\Phi(X)$ & $:=$ & $\left\{
\begin{array}
[c]{l}%
\mathbb{E}%
{\displaystyle\int\nolimits_{0}^{T}}
\varphi(X_{t})dt,\,\,\text{\textit{if }}\varphi\left(  X\right)  \in
L^{1}(\Omega\times]0,T[),\smallskip\vspace{0.04in}\\
+\infty,\;\text{\textit{otherwise}}%
\end{array}
\right.  $\vspace{0.04in}$\vspace{0.04in}$%
\end{tabular}
\ \ \ \ \ \ \ \ \ \ \ \ \ \
\]
For each%
\[
\left(  U,V\right)  \in\mathbb{D}:=Dom\left(  \Phi\right)  \times L^{2}\left(
\Omega\times\left[  0,T\right]  ;\mathbb{R}^{d}\right)
\]
we introduce the function%
\[
J_{\left(  U,V\right)  }:\mathbb{S}:=L^{2}(\Omega,\mathcal{F}_{T}%
,\mathbb{P},\mathbb{R}^{d})\times\Lambda_{d}^{2}\left(  0,T\right)  \times
S_{d}^{2}\left(  0,T\right)  \times\Lambda_{d\times k}^{2}\left(  0,T\right)
\rightarrow\mathbb{R}%
\]
by%
\[%
\begin{array}
[c]{r}%
J_{\left(  U,V\right)  }(\eta,G,Y,Z):=\dfrac{1}{2}\mathbb{E}\left\vert
\eta-\xi\right\vert ^{2}+\mathbb{E}%
{\displaystyle\int_{0}^{T}}
\left\langle U_{t}-Y_{t},F(t,U_{t},V_{t})-G_{t}\right\rangle dt\medskip\\
-\dfrac{1}{2}\mathbb{E}%
{\displaystyle\int_{0}^{T}}
\left\vert Z_{t}-V_{t}\right\vert ^{2}dt+\Phi\left(  Y\right)  -\Phi\left(
U\right)
\end{array}
\]
and consider the functional $\hat{J}:\mathbb{S}\rightarrow]-\infty,+\infty],$%
\[
\hat{J}(\eta,G,Y,Z):=\sup_{\left(  U,V\right)  \in\mathbb{D}}J_{\left(
U,V\right)  }(\eta,G,Y,Z).
\]

\begin{remark}
$\hat{J}:\mathbb{S}\rightarrow]-\infty,+\infty]$ is a lower semicontinuous
function as supremum of lower semicontinuous functions.
\end{remark}

We now have

\begin{proposition}
The mapping $\hat{J}$ has the following properties:

\begin{itemize}
\item[$\left(  a\right)  $] $\ \hat{J}(\eta,G,Y,Z)\geq0,\ \forall
(\eta,G,Y,Z)\in\mathbb{S}$ and $\hat{J}$ is not identical $+\infty.$

\item[$\left(  b\right)  $] $\ $Let $(\hat{\eta},\hat{G},\hat{Y},\hat{Z}%
)\in\mathbb{S}.$ Then
\[
\hat{J}(\hat{\eta},\hat{G},\hat{Y},\hat{Z})=0\quad\text{iff}\quad\hat{\eta
}=\xi,\;F(\hat{Y},\hat{Z})-\hat{G}\in\partial\Phi(\hat{Y}).
\]

\item[$\left(  c\right)  $] $\ $The restriction of $\hat{J}$ to the closed
convex set%
\[
\mathbb{K}=\left\{  (\eta,G,Y,Z)\in\mathbb{S}:\;Y_{t}=\eta+%
{\displaystyle\int_{t}^{T}}
G_{s}ds-%
{\displaystyle\int_{t}^{T}}
Z_{s}dB_{s},\ \forall t\in\left[  0,T\right]  \right\}
\]
is a convex lower semicontinuous function. If $(\hat{\eta},\hat{G},\hat
{Y},\hat{Z})\in\mathbb{K},$ then%
\[
\hat{J}(\hat{\eta},\hat{G},\hat{Y},\hat{Z})=0\quad\text{iff}\quad\hat{\eta
}=\xi\;\text{and }(\hat{Y},\hat{Z},\hat{H}),\text{with }\hat{H}=F(\hat{Y}%
,\hat{Z})-\hat{G}%
\]
is a solution of the BSVI (\ref{bsvi1}).
\end{itemize}
\end{proposition}

\begin{proof}
$\left(  a\right)  \;$If $Y\notin Dom\left(  \Phi\right)  $ then $J_{\left(
U,V\right)  }(\eta,G,Y,Z)=+\infty$ and if $Y\in Dom\left(  \Phi\right)  $, we
have $\hat{J}(\eta,G,Y,Z)\geq J_{(Y,Z)}(\eta,G,Y,Z)\geq0$. Moreover, $\hat{J}$
is a proper function since for $v_{0}\in\partial\varphi\left(  u_{0}\right)  $
and $\eta^{0}=\xi,$ $Y_{t}^{0}=u_{0},$ $Z_{t}^{0}=0$, $G_{t}^{0}=F\left(
t,u_{0},0\right)  -v_{0}$ we have (using the assumption (\ref{bsvi-a1})) that%
\[
\hat{J}_{\left(  U,V\right)  }\left(  \eta^{0},G^{0},Y^{0},Z^{0}\right)
\leq0,\;\text{for all }\left(  U,V\right)  \in\mathbb{D}.
\]

$\left(  b\right)  \;$If $\hat{J}(\hat{\eta},\hat{G},\hat{Y},\hat{Z})=0$%
\ then%
\[
J_{\left(  U,V\right)  }(\hat{\eta},\hat{G},\hat{Y},\hat{Z})\leq0,\ \forall
U\in Dom\left(  \Phi\right)  ,\ \forall V\in L^{2}\left(  \Omega\times\left[
0,T\right]  ;\mathbb{R}^{d}\right)  .
\]
So, for all $\left(  U,V\right)  \in\mathbb{D}$,%
\[
\dfrac{1}{2}\mathbb{E}\left\vert \hat{\eta}-\xi\right\vert ^{2}+\mathbb{E}%
\int_{0}^{T}\left\langle U_{t}-\hat{Y}_{t},F(U_{t},V_{t})-\hat{G}%
_{t}\right\rangle dt-\dfrac{1}{2}\mathbb{E}\int_{0}^{T}|\hat{Z}_{t}-V_{t}%
|^{2}dt+\Phi(\hat{Y})-\Phi\left(  U\right)  \leq0,
\]
which yields $\hat{Y}\in Dom\left(  \Phi\right)  $; taking in particular
$U=\hat{Y}$ and $V=\hat{Z}$, we infer%
\[
\hat{\eta}=\xi,\ \mathbb{\mathbb{P}}\text{-}a.s.
\]
Hence, for all $\left(  U,V\right)  \in\mathbb{D}$,%
\begin{equation}
\mathbb{E}\int_{0}^{T}\left\langle U_{t}-\hat{Y}_{t},F(U_{t},V_{t})-\hat
{G}_{t}\right\rangle dt+\Phi(\hat{Y})\leq\dfrac{1}{2}\mathbb{E}\int_{0}%
^{T}|\hat{Z}_{t}-V_{t}|^{2}dt+\Phi\left(  U\right)  . \label{inegidentif}%
\end{equation}
Since $\mathbb{D}$ is a convex set, we can replace $\left(  U,V\right)  $ by
$(\left(  1-\lambda\right)  \hat{Y}+\lambda U,\left(  1-\lambda\right)
\hat{Z}+\lambda V),$ where $\lambda\in\left(  0,1\right)  .$ The convexity of
$\Phi$ leads to the following inequality%
\[%
\begin{array}
[c]{r}%
\mathbb{E}%
{\displaystyle\int_{0}^{T}}
\left\langle U_{t}-\hat{Y}_{t},F(\left(  1-\lambda\right)  \hat{Y}_{t}+\lambda
U_{t},\left(  1-\lambda\right)  \hat{Z}_{t}+\lambda V_{t})-\hat{G}%
_{t}\right\rangle dt\medskip\\
\leq\dfrac{\lambda}{2}\mathbb{E}%
{\displaystyle\int_{0}^{T}}
|\hat{Z}_{t}-V_{t}|^{2}dt+\Phi\left(  U\right)  -\Phi(\hat{Y}).
\end{array}
\]
Passing to $\liminf\limits_{\lambda\rightarrow0}$, we deduce%
\[
\mathbb{E}%
{\displaystyle\int_{0}^{T}}
\left\langle U_{t}-\hat{Y}_{t},F(\hat{Y}_{t},\hat{Z}_{t})-\hat{G}%
_{t}\right\rangle dt+\Phi(\hat{Y})\leq\Phi\left(  U\right)  ,\ \forall U\in
Dom\left(  \Phi\right)  ,
\]
that is%
\[
F(\hat{Y},\hat{Z})-\hat{G}\in\partial\Phi(\hat{Y}).
\]
Conversely, using assumption (\ref{bsvi-a1}) we have%
\begin{align*}
&  J_{\left(  U,V\right)  }(\xi,\hat{G},\hat{Y},\hat{Z})\\
&  =\mathbb{E}\int_{0}^{T}\left\langle U_{t}-\hat{Y}_{t},F(U_{t},V_{t}%
)-\hat{G}_{t}\right\rangle dt-\dfrac{1}{2}\mathbb{E}\int_{0}^{T}|\hat{Z}%
_{t}-V_{t}|^{2}dt+\Phi(\hat{Y})-\Phi\left(  U\right) \\
&  \leq\mathbb{E}\int_{0}^{T}\left\langle U_{t}-\hat{Y}_{t},F(U_{t}%
,V_{t})-F(\hat{Y}_{t},\hat{Z}_{t})\right\rangle dt-\dfrac{1}{2}\mathbb{E}%
\int_{0}^{T}|\hat{Z}_{t}-V_{t}|^{2}dt\\
&  +\mathbb{E}\int_{0}^{T}\left\langle U_{t}-\hat{Y}_{t},F(\hat{Y}_{t},\hat
{Z}_{t})-\hat{G}_{t}\right\rangle dt+\Phi(\hat{Y})-\Phi\left(  U\right) \\
&  \leq0
\end{align*}
and, consequently, $\hat{J}(\hat{\eta},\hat{G},\hat{Y},\hat{Z})=0.\smallskip$

$\left(  c\right)  \;$If, moreover, $(\hat{\eta},\hat{G},\hat{Y},\hat{Z}%
)\in\mathbb{K},$ then%
\[
Y_{t}+%
{\displaystyle\int_{t}^{T}}
(F(\hat{Y}_{s},\hat{Z}_{s})-\hat{G}_{s})ds=\hat{\eta}+%
{\displaystyle\int_{t}^{T}}
F(\hat{Y}_{s},\hat{Z}_{s})ds-%
{\displaystyle\int_{t}^{T}}
Z_{s}dB_{s}%
\]
and%
\[
F(\hat{Y},\hat{Z})-\hat{G}\in\partial\Phi(\hat{Y}),
\]
that is, $(\hat{Y},\hat{Z},F(\hat{Y},\hat{Z})-\hat{G})$ is solution of the SVI
(\ref{bsvi1}).

The convexity of $\hat{J}$ on $\mathbb{K}$ is obtained as follows: by the
Energy Equality we have%
\[
\left\vert Y_{0}\right\vert ^{2}+\mathbb{E}\int_{0}^{T}\left\vert
Z_{s}\right\vert ^{2}ds=\mathbb{E}\left\vert \eta\right\vert ^{2}%
+2\mathbb{E}\int_{0}^{T}\left\langle Y_{s},G_{s}\right\rangle ds
\]
and $J_{(U,V)}(\eta,G,Y,Z)$ becomes%
\begin{align*}
J_{(U,V)}(\eta,G,Y,Z)  &  =\dfrac{1}{2}\mathbb{E}\left\vert \eta
-\xi\right\vert ^{2}+\mathbb{E}%
{\displaystyle\int_{0}^{T}}
\left\langle U_{t}-Y_{t},F(U_{t},V_{t})\right\rangle dt-\mathbb{E}%
{\displaystyle\int_{0}^{T}}
\left\langle U_{t},G_{t}\right\rangle dt\\
&  +\mathbb{E}%
{\displaystyle\int_{0}^{T}}
\left\langle Y_{t},G_{t}\right\rangle dt-\dfrac{1}{2}\mathbb{E}\int_{0}%
^{T}\left\vert Z_{t}-V_{t}\right\vert ^{2}dt+\Phi\left(  Y\right)
-\Phi\left(  U\right) \\
&  =\dfrac{1}{2}\mathbb{E}\left\vert \xi\right\vert ^{2}-\mathbb{E}%
\left\langle \eta,\xi\right\rangle +\mathbb{E}\int_{0}^{T}\left\langle
U_{t}-Y_{t},F(U_{t},V_{t})\right\rangle dt-\mathbb{E}\int_{0}^{T}\left\langle
U_{t},G_{t}\right\rangle dt\\
&  +\mathbb{E}%
{\displaystyle\int_{0}^{T}}
\left\langle Z_{t},V_{t}\right\rangle dt-\dfrac{1}{2}\mathbb{E}%
{\displaystyle\int_{0}^{T}}
\left\vert V_{t}\right\vert ^{2}dt+\dfrac{1}{2}\mathbb{E}\left\vert
Y_{0}\right\vert ^{2}+\Phi\left(  Y\right)  -\Phi\left(  U\right)  .
\end{align*}
Hence $\hat{J}$ is a convex l.s.c. function as supremum of convex l.s.c. functions.

The proof is now complete.\hfill
\end{proof}

\begin{acknowledgement}
The authors are grateful to the referees for the attention in reading this
paper and for their very useful suggestions.
\end{acknowledgement}

\end{document}